\documentclass[a4paper,10pt]{article}

\usepackage{inputenc, amsmath, amsthm, amssymb}
\usepackage{latexsym}

\usepackage{mathrsfs}

\usepackage{a4wide}

\usepackage{fancyhdr, color, multicol, graphicx}

\numberwithin{equation}{section}

\newtheorem{thm}{Theorem}[section]
\newtheorem{lem}[thm]{Lemma}

\newtheorem{cor}[thm]{Corollary}
\newtheorem{defi}[thm]{Definition}

\theoremstyle{definition}

\newenvironment{prf}{\textit{Proof}.}
{\hfill $\Box$}

\newcommand{\eps}{\varepsilon}

\begin{document}

\begin{center}
\Large \textbf{Locally Supported Wavelets for the Separation of Spherical Vector Fields with Respect to their Sources}
\\[3ex]\normalsize C. Gerhards \footnotetext{\begin{flushleft}University of Kaiserslautern, Geomathematics Group, 67653 Kaiserslautern, PO Box 3049, Germany
\\e-mail: gerhards@mathematik.uni-kl.de\end{flushleft}}
\end{center}

\textbf{Abstract} We provide a space domain oriented separation of magnetic fields into parts generated by sources in the exterior and sources in the interior of a given sphere. The separation itself is well-known in geomagnetic modeling, usually in terms of a spherical harmonic analysis or a wavelet analysis that is spherical harmonic based. In contrast to these frequency oriented methods, we use a more spatially oriented approach in this paper. We derive integral representations with explicitly known convolution kernels. Regularizing these singular kernels allows a multiscale representation of the internal and external contributions to the magnetic field with locally supported wavelets. This representation is applied to a set of CHAMP data for crustal field modeling.
\\

\\


\section{Introduction}

The Earth's magnetic field is a complex structure consisting of various contributions, such as the dominating core field, the crustal field, and effects from iono- and magnetospheric processes. A major task in understanding the geomagnetic field is the separation of these contributions. An overview on different approaches to this is given, e.g., in \cite{olsen10}. A first step is the mathematical separation of magnetic field measurements taken at satellite altitude into contributions from sources in the exterior of the orbit and contributions from sources in the interior. Generally, we assume the magnetic field $b$ to be governed by the pre-Maxwell equations
\begin{eqnarray*}
\nabla\wedge b&=&\mu_0j,
\\\nabla\cdot b&=&0,
\end{eqnarray*}
with $j$ describing the source current density and $\mu_0$ the vacuum permeability ($\wedge$ denotes the vector product). If no source currents $j$ are present, one has $b=\nabla U$, for some harmonic potential $U$, and the typical approach to modeling the magnetic field is the so-called Gauss representation of the corresponding potential in terms of scalar spherical harmonics $Y_{n,k}$ (see, e.g., \cite{back96} and \cite{sabaka10}). Generally, however, satellite data is collected in a source region of the magnetic field. Then the Mie decomposition allows a decomposition of the magnetic field into a poloidal part $p_b$ and a toroidal part $q_b$. The toroidal part describes the magnetic field due to poloidal current densities $p_j$, while the poloidal part can be split into a part $p_b^{ext}$ that is due toroidal sources in the exterior of the satellite's orbit, and a part $p_b^{int}$ that is due to toroidal sources in the interior. A more detailed description can be found, e.g., in \cite{back86} and \cite{back96}. In this setting, the quantities $p_b^{int}$, $p_b^{ext}$ and $q_b$ can be expanded in a system of vector spherical harmonics $\tilde{y}_{n,k}^{(1)}$, $\tilde{y}_{n,k}^{(2)}$ and $\tilde{y}_{n,k}^{(3)}$, respectively (a system that actually originates in quantum mechanics; see, e.g., \cite{edmonds57}). 

However, due to the global nature of scalar and vector spherical harmonics, they are not the best choice for modeling strongly localized structures, such as the Earth's crustal field, or modeling from only locally available data. Several multiscale approaches with spatially better localizing kernels have been developed to improve this drawback, e.g., in \cite{cham05}, \cite{hol03} for potential fields, in \cite{may06}, \cite{maymai} for the above described separation with respect to the sources, and in \cite{bay01}, \cite{mai05} for a representation of ionospheric magnetic fields and current densities. A comprehensive introduction of kernel functions for such methods can also be found in \cite{free98}. 

It is the aim of this paper to transfer the multiscale approach described in \cite{may06}, \cite{maymai}, which is based on a construction of scaling and wavelet kernels in frequency domain (i.e., based on an adequate superposition of the vector spherical harmonics $\tilde{y}_{n,k}^{(i)}$), to a setting where the scaling and wavelet kernels are constructed entirely in space domain. For that purpose, the vector spherical harmonics $\tilde{y}_{n,k}^{(i)}$, $i=1,2,3$, are described by operators $\tilde{o}^{(i)}$, $i=1,2,3$. A decomposition of the magnetic field in terms of these operators, in combination with the spherical Helmholtz decomposition, allows an integral expression of the quantities $p_b^{int}$, $p_b^{ext}$, $q_b$. Motivated by \cite{free06}, a regularization of the convolution kernels appearing in this integral expression provides a multiscale representation with wavelets that are locally supported in space. The multiscale representation is described in detail in Section \ref{sec:multiscale}. There, we also apply the derived algorithm to a set of real CHAMP satellite data. The preparatory construction of the regularized kernels and a decomposition with respect to the operators $\tilde{o}^{(i)}$ is described in Sections \ref{sec:kernels} and \ref{sec:decomps}. Section \ref{sec:basics} provides fundamental aspects on Legendre polynomials and scalar and vector spherical harmonics.

\section{Preliminaries}\label{sec:basics}

By $P_n:[-1,1]\to\mathbb{R}^3$, $n\in\mathbb{N}_0$, we denote the set of Legendre polynomials of degree $n$, by $Y_{n,k}:\Omega\to\mathbb{R}$, $n\in\mathbb{N}_0$, $k=1,\ldots,2n+1$, an orthonormal set of spherical harmonics of degree $n$ and order $k$ ($\Omega_R=\{x\in\mathbb{R}^3|\,|x|=R\}$ denotes the sphere of radius $R>0$ and $\Omega=\Omega_1$ the unit sphere). The fundamental connection between these two function systems is the so-called addition theorem,
\begin{eqnarray*}
\sum_{k=1}^{2n+1}Y_{n,k}(\xi)Y_{n,k}(\eta)=\frac{2n+1}{4\pi}P_n(\xi\cdot\eta),\quad\xi,\eta\in\Omega.
\end{eqnarray*}
This allows us to expand zonal kernels (i.e., functions $F:\Omega\times\Omega\to\mathbb{R}$ that satisfy  $F(\xi,\eta)=G(\xi\cdot\eta)$, $\xi,\eta\in\Omega$, for an adequate function $G:[-1,1]\to\mathbb{R}$) in terms of Legendre polynomials. Known closed representations for certain series of Legendre polynomials can then be used to derive closed representations for some zonal kernels appearing in this paper. One of these series is the generating series for the Legendre polynomials,
\begin{eqnarray*}
\sum_{n=0}^\infty h^nP_n(t)=\frac{1}{\sqrt{1+h^2-2ht}},\quad t\in[-1,1],\,h\in(-1,1).
\end{eqnarray*}
From this, one can derive various further representations that are, e.g., listed in \cite{hansen}. Of importance to us are the following ones.

\begin{lem}\label{lem:genseriesgreenfunc}
For $t\in(-1,1)$, we have
\begin{eqnarray*}
\sum_{n=1}^\infty \frac{1}{n}P_n(t)&=&\ln\left(\frac{\sqrt{2}\sqrt{1-t}-1+t}{1-t^2}\right)+\ln\left(2\right),\label{eqn:genseriesgreenfunc1hto0}
\\\sum_{n=1}^\infty \frac{1}{n+1}P_n(t)&=&\ln\left(1+\frac{\sqrt{2}}{\sqrt{1-t}}\right)-1.\label{eqn:genseriesgreenfunc2hto0}
\end{eqnarray*}
\end{lem}

Furthermore, the generating series for the Legendre polynomials yields an expansion of the single layer kernel. This is of interest since it allows an integral definition of the single layer operator and a definition in terms of pseudodifferential operators.

\begin{lem}\label{lem:seriespotfunc}
Let $x,y\in\mathbb{R}^3$ with $|x|<|y|$. Then
\begin{equation*}
\frac{1}{|x-y|}=\frac{1}{|y|}\sum_{n=0}^\infty\left(\frac{|x|}{|y|}\right)^nP_n\left(\frac{x}{|x|}\cdot\frac{y}{|y|}\right).
\end{equation*}
\end{lem}

The set of the previously mentioned spherical harmonics yields a complete orthonormal system in $L^2(\Omega)=\{F:\Omega\to\mathbb{R}|\int_\Omega |F(\eta)|^2d\omega(\eta)<\infty\}$. The modeling of magnetic fields, however, is in first place a vectorial problem. For that purpose, we introduce two different complete sets of vector spherical harmonics. The first set requires the operators
\begin{eqnarray}
o^{(1)}_\xi F(\xi)&=&\xi F(\xi),\quad\xi\in\Omega,\label{eqn:o1operator}
\\o^{(2)}_\xi F(\xi)&=&\nabla^*_\xi F(\xi),\quad\xi\in\Omega,\label{eqn:o2operator}
\\o^{(3)}_\xi F(\xi)&=&L^*_\xi F(\xi),\quad\xi\in\Omega,\label{eqn:o3operator}
\end{eqnarray}  
for $F:\Omega\to\mathbb{R}$ a sufficiently smooth scalar function, $\nabla^*$ the surface gradient (i.e., the tangential part of the gradient $\nabla$; more precisely, $\nabla_x=\xi\frac{\partial}{\partial r}+\frac{1}{r}\nabla_\xi^*$, for $x=r\xi\in\mathbb{R}^3$, with $r=|x|$, $\xi=\frac{x}{|x|}$), and $L^*$ the surface curl gradient (acting as $L_\xi^*=\xi\wedge\nabla_\xi^*$ in a point $\xi\in\Omega$). A complete orthonormal system in $l^2(\Omega)=\{f:\Omega\to\mathbb{R}^3|$ $\int_\Omega|f(\eta)|^2d\omega(\eta)<\infty\}$ is then given via
\begin{eqnarray}
y_{n,k}^{(i)}=(\mu_n^{(i)})^{-\frac{1}{2}}o^{(i)}Y_{n,k},\quad i=1,2,3,\,n\in\mathbb{N}_{0_i},\,k=1,\ldots,2n+1,
\end{eqnarray}
where $0_i$ is an abbreviation for $0_1=0$ and $0_i=1$, $i=2,3$, and $\mu_n^{(i)}$ denotes the normalization constants $\mu_n^{(1)}=1$ and $\mu_n^{(i)}=n(n+1)$, $i=2,3$. Concerning the notation, upper case letters, such as $F$, $Y_{n,k}$, generally denote scalar valued functions, lower case letters, such as $f$, $y_{n,k}^{(i)}$, denote vector valued functions, and bold face letters denote tensor valued functions. The same notation holds for the function spaces $C^{(k)}(\Omega)$, $c^{(k)}(\Omega)$ of $k$-times continuously differentiable functions and the spaces $L^2(\Omega)$, $l^2(\Omega)$ of square integrable functions.

The second set of vector spherical harmonics requires the modified operators 
\begin{eqnarray}
\tilde{o}^{(1)}&=&o^{(1)}\left(D+\textnormal{\footnotesize $\frac{1}{2}$}\right)-o^{(2)},\label{eqn:to1operator}
\\\tilde{o}^{(2)}&=&o^{(1)}\left(D-\textnormal{\footnotesize $\frac{1}{2}$}\right)+o^{(2)},\label{eqn:to2operator}
\\\tilde{o}^{(3)}&=&o^{(3)},\label{eqn:to3operator}
\end{eqnarray}
where
\begin{eqnarray}
D=\left(-\Delta^*+\textnormal{\footnotesize $\frac{1}{4}$}\right)^{\frac{1}{2}}.\label{eqn:d}
\end{eqnarray}
By $\Delta^*$ we denote the Beltrami operator $\nabla^*\cdot\nabla^*$. The operator $D$ is treated in more detail in Subsection \ref{sec:singlelayer}. A second complete orthonormal system in $l^2(\Omega)$ is then given via
\begin{eqnarray}
\tilde{y}_{n,k}^{(i)}=(\tilde{\mu}_n^{(i)})^{-\frac{1}{2}}\tilde{o}^{(i)}Y_{n,k},\quad i=1,2,3,\,n\in\mathbb{N}_{0_i},\,k=1,\ldots,2n+1,\label{eqn:ytilde}
\end{eqnarray}
where $\tilde{\mu}_n^{(i)}$ denotes the normalization constants $\tilde{\mu}_n^{(1)}=(n+1)(2n+1)$, $\tilde{\mu}_n^{(2)}=n(2n+1)$ and $\tilde{\mu}_n^{(3)}=n(n+1)$. The advantage of this basis system is its connection to the inner and outer harmonics, i.e., the functions $H_{n,k}^{int}(x)=\frac{1}{R}\big(\frac{|x|}{R}\big)^nY_{n,k}\big(\frac{x}{|x|}\big)$, $x\in{\Omega_R^{int}}=\{x\in\mathbb{R}^3|\,|x|<R\}$, and $H_{n,k}^{ext}(x)=\frac{1}{R}\big(\frac{R}{|x|}\big)^{n+1}Y_{n,k}\big(\frac{x}{|x|}\big)$, $x\in{\Omega^{ext}_R}=\{x\in\mathbb{R}^3|\,|x|>R\}$,  which yield solutions to the inner and outer Dirichlet boundary value problem, respectively (i.e., boundary values $H_{n,k}^{int}=H_{n,k}^{ext}=Y_{n,k}$ on $\Omega_R$ and $\Delta H_{n,k}^{int}=0$ in $\Omega_R^{int}$, $\Delta H_{n,k}^{ext}=0$ in $\Omega_R^{ext}$). We have
\begin{eqnarray}
\nabla_xH^{int}_{n,k}(x)&=&\frac{1}{R^2}\left(\frac{r}{R}\right)^{n-1}(\tilde{\mu}_n^{(2)})^{\frac{1}{2}}\tilde{y}_{n,k}^{(2)}(\xi),\qquad r=|x|,\,x=r\xi\in\overline{\Omega_R^{int}},\label{eqn:hint}
\\-\nabla_xH^{ext}_{n,k}(x)&=&\frac{1}{R^2}\left(\frac{R}{r}\right)^{n+2}(\tilde{\mu}_n^{(1)})^{\frac{1}{2}}\tilde{y}_{n,k}^{(1)}(\xi),\qquad r=|x|,\, x=r\xi\in\overline{\Omega_R^{ext}}.\label{eqn:hext}
\end{eqnarray}
For a more comprehensive introduction of the function systems mentioned in this section, the reader is referred to, e.g., \cite{free98} and the references therein. The special importance of the last set of vector spherical harmonics in geomagnetic modeling is well emphasized, e.g., in \cite{back96}, \cite{may06} and \cite{maymai}. In this paper, however, they are only to be understood as a motivation for the Helmholtz decomposition and a modified decomposition with respect to $\tilde{o}^{(i)}$. Our main goal is to actually avoid spherical harmonic representations.

\section{Regularized Kernels}\label{sec:kernels}

Green's function for the Beltrami operator and the single layer kernel are especially useful when working with differential equations involving the operators $\nabla^*$, $L^*$, $\Delta^*$ and $D$. We briefly recapitulate some of the properties of these functions and the corresponding operators before we introduce a regularization for both kernels separately and for their combination. To achieve integral representations for the scalars of the classical Helmholtz decomposition, it is actually sufficient to only have Green's function. The single layer kernel becomes necessary when we introduce a decomposition that pays tribute to interior and exterior sources.

\subsection{Green's Function}\label{sec:greenfunc}

By \emph{Green's function with respect to the Beltrami operator} we denote the uniquely defined function $G(\Delta^*;\cdot):[-1,1)\to\mathbb{R}$ satisfying the properties
\begin{itemize}
\item[(i)] $\eta\mapsto G(\Delta^*;\xi\cdot\eta)$ is twice continuously differentiable on the set $\{\eta\in\Omega|$ $1-\xi\cdot\eta>0\}$, and
\vspace*{-1mm}
\[
\Delta^*_\eta G(\Delta^*;\xi\cdot\eta)= -\frac{1}{4 \pi}, \quad 1-\xi\cdot\eta>0,
\]
for any fixed $\xi \in \Omega$,
\item[(ii)]for any fixed $\xi \in \Omega$, the function
\vspace*{-1mm}
\[
\eta\mapsto G(\Delta^*;\xi\cdot\eta)- \frac{1}{4 \pi}\ln(1-\xi \cdot \eta),
\]
is continuously differentiable on $\Omega$,
\item[(iii)] for any fixed $\xi \in \Omega$,
\vspace*{-1mm}
\[
\frac{1}{4\pi}\int_\Omega G(\Delta^*;\xi\cdot\eta)d\omega(\eta) = 0 .
\]
\end{itemize}

One can verify the following explicit representation,
\begin{eqnarray}\label{eqn:beltramigreenfct}
G(\Delta^*;\xi\cdot\eta)&=&\frac{1}{4\pi}\ln(1-\xi\cdot\eta)+\frac{1}{4\pi}(1-\ln(2)),\quad1-\xi\cdot\eta>0.\label{eqn:repmodgreenfunc}
\end{eqnarray}
The bilinear series expansion reads
\begin{eqnarray}
G(\Delta^*; \xi\cdot\eta)&=&\sum_{n=1}^\infty\sum_{k=1}^{2n+1}\frac{1}{-n(n+1)}Y_{n,k}(\xi)Y_{n,k}(\eta),\quad 1-\xi\cdot\eta>0.
\end{eqnarray}
Observing that $\eta\mapsto\Delta^*_\eta G(\Delta^*;\xi\cdot\eta)$ only varies by the constant $-\frac{1}{4 \pi}$ from the Dirac distribution motivates the following theorems which express a sufficiently smooth function by its integral mean value and a correction term involving Green's function. For more details, the reader is again referred to \cite{free98} and the references therein.

\begin{thm}[Fundamental Theorem for $\Delta^*$]\label{thm:fundthmdelta}
Let $F$ be of class $C^{(2)}(\Omega)$. Then
\begin{equation*}
F(\xi)=\frac{1}{4\pi}\int_\Omega F(\eta)d\omega(\eta)+\int_\Omega G(\Delta^*;\xi\cdot\eta)\Delta^*_\eta F(\eta)\,d\omega(\eta),\quad\xi\in\Omega.
\end{equation*}
\end{thm}

\begin{thm}[Fundamental Theorem for $\nabla^*$ and $L^*$]\label{thm:fundthmgradcurl}
Let $F$ be of class $C^{(1)}(\Omega)$. Then
\begin{eqnarray*}
F(\xi)&=&\frac{1}{4\pi}\int_\Omega F(\eta)d\omega(\eta)-\int_\Omega \Lambda_\eta^*G(\Delta^*;\xi\cdot\eta)\cdot \Lambda_\eta^*F(\eta)\,d\omega(\eta),\quad\xi\in\Omega,
\end{eqnarray*}
where $\Lambda^*$ denotes one of the operators $\nabla^*$ or $L^*$.
\end{thm}

These theorems directly yield simple integral representations for solutions to the spherical differential equations with respect to $\nabla^*$, $L^*$, and $\Delta^*$.

Next, we present a spatial regularization of $G(\Delta^*;\cdot)$ around its singularity. This is a crucial step for the later definition of the scaling and wavelet kernels of the multiscale representation.

\begin{defi}[Regularized Green's Function]\label{def:reggreen}
Let $R^\rho$, $\rho>0$, be of class $C^{(n)}([-1,1])$, $n\in\mathbb{N}$ fixed, satisfying
\begin{equation*}
\lim_{\rho\to0+}\rho^{\frac{k}{2}}\int_{1-\rho}^1\left|\left(\frac{d}{dt}\right)^kR^\rho(t)\right|\,dt=0,\quad k=0,1,\label{eqn:limitrrho}
\end{equation*}
and 
\begin{equation*}
\left[\left(\frac{d}{dt}\right)^k R^\rho(t)\right]_{t=1-\rho}=\left[\left(\frac{d}{dt}\right)^kG(\Delta^*;t)\right]_{t=1-\rho},\quad k=0,1,\ldots,n. \end{equation*}
Then the function
\begin{eqnarray*}
G^\rho(\Delta^*;\xi\cdot\eta)=\left\{\begin{array}{ll}
G(\Delta^*;\xi\cdot\eta),&1-\xi\cdot\eta\geq\rho,
\\[1.25ex]R^\rho(\xi\cdot\eta),&1-\xi\cdot\eta<\rho,
\end{array}\right.
\end{eqnarray*}
is called \emph{regularized Green's function (of order $n$)}. $R^\rho$ is called the regularization function.
\end{defi}

\begin{figure}
\begin{center}
\scalebox{0.37}{\includegraphics{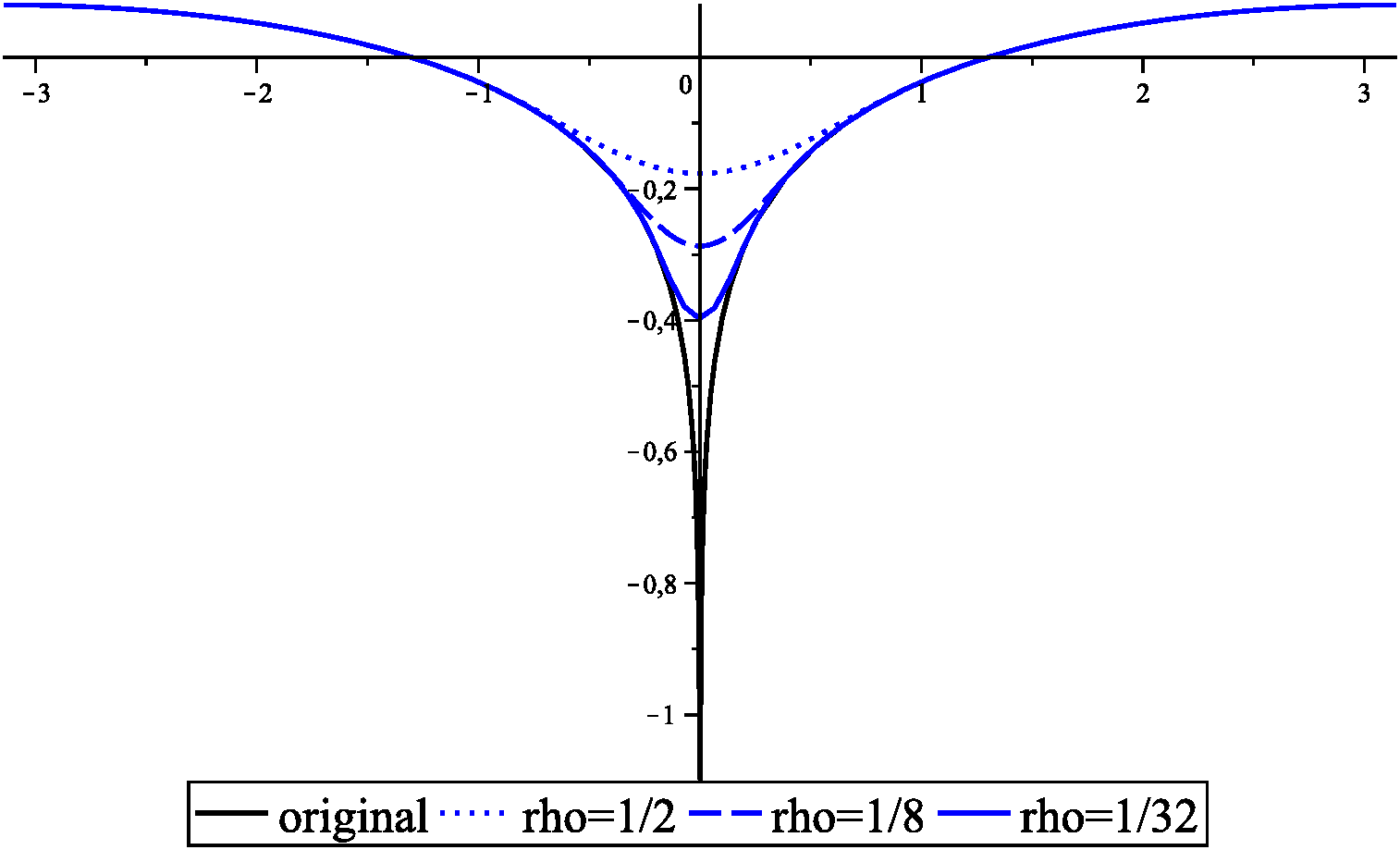}}\qquad \scalebox{0.36}{\includegraphics{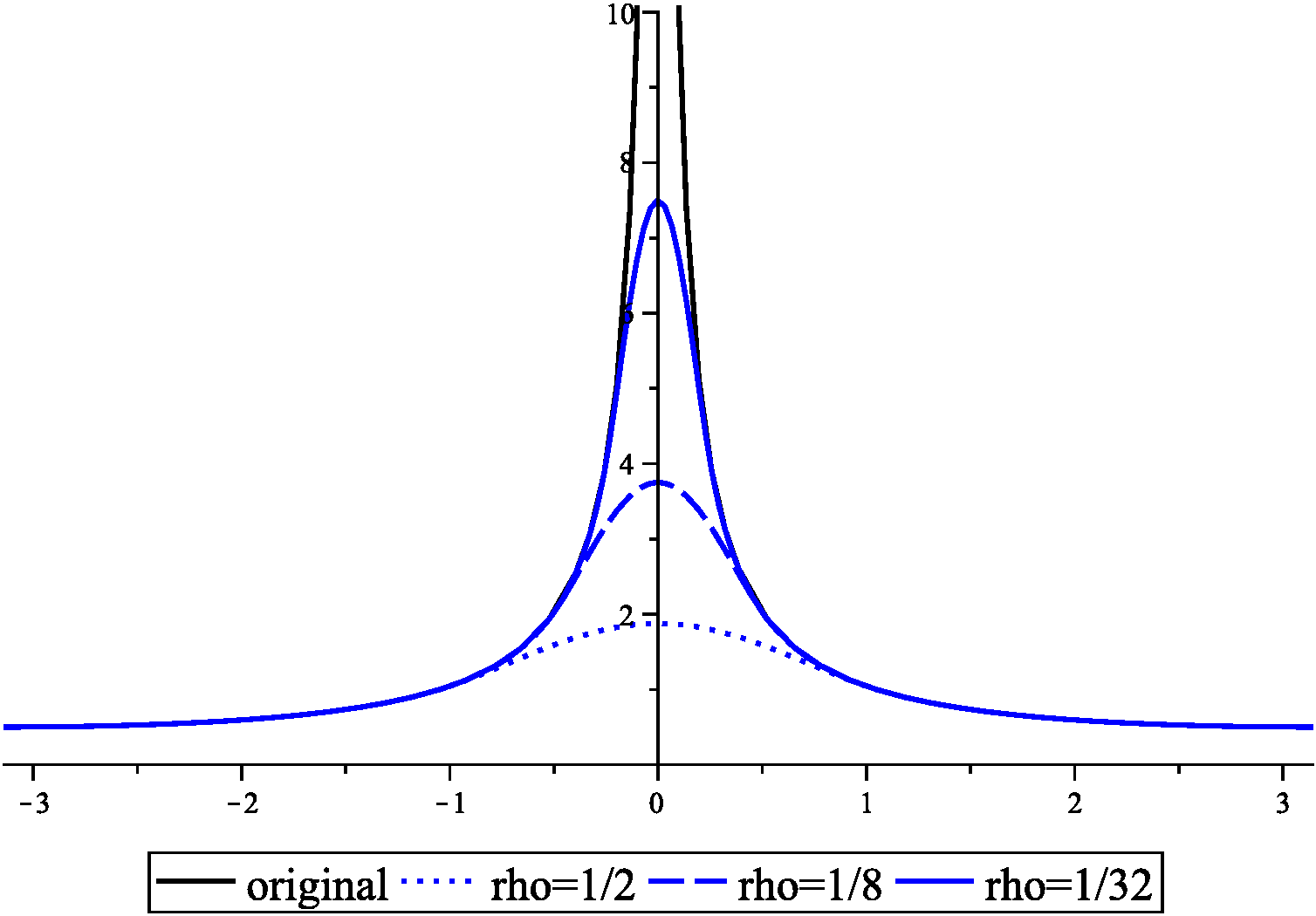}}
\end{center}
\caption{Plot of a twice continuously differentiable regularization $\theta\mapsto G^\rho(\Delta^*;\cos(\theta))$ (left) and a twice continuously differentiable regularization $\theta\mapsto S^\rho(\cos(\theta))$ (right) at different scales $\rho$.}\label{fig:reggreenfunc}
\end{figure}

A typical choice for $R^\rho$ is the Taylor series of $G(\Delta^*;\cdot)$ centered at $1-\rho$ and truncated at some power $n$. An exemplary plot for different scaling parameters $\rho$ can be found in Figure \ref{fig:reggreenfunc}. Similar regularizations, but only for Taylor polynomials up to degree $2$, have been used in other areas of geosciences, e.g., in \cite{fehl07}, \cite{fehl08}, and \cite{free06}. To be able to state a multiscale decomposition, it has to be guaranteed that convolutions with the regularized kernels converge to convolutions with the original kernels. The proofs are based on the fact that $\eta\mapsto G(\Delta^*;\xi\cdot\eta)$, $\eta\mapsto \nabla_\xi^*G(\Delta^*;\xi\cdot\eta)$ and $\eta\mapsto L_\xi^*G(\Delta^*;\xi\cdot\eta)$ are integrable on the sphere $\Omega$, uniformly with respect to $\xi\in\Omega$, and can be found in \cite{free09} and \cite{freeger10}.

\begin{lem}\label{lem:convg}
Let  $G^\rho(\Delta^*; \cdot)$ be of class $ C^{(1)}([-1, 1])$ and $F$ of class $C^{(0)}(\Omega)$. Then we have
\begin{eqnarray*}
\lim_{\rho\to0+}\sup_{\xi\in\Omega}\left|\int_\Omega G^\rho(\Delta^*;\xi\cdot\eta)F(\eta)d\omega(\eta)-\int_\Omega G(\Delta^*;\xi\cdot\eta)F(\eta)d\omega(\eta)\right|=0.
\end{eqnarray*}
\end{lem}

\begin{lem}\label{lem:convdiffg}
Let $F$ be of class $C^{(0)}(\Omega)$ and $G^\rho(\Delta^*; \cdot)$ of class $ C^{(1)}([-1, 1])$. Then
\begin{eqnarray*}
\lim_{\rho\to0+}\sup_{\xi\in\Omega}\left|\int_\Omega \Lambda^*_\xi G^\rho(\Delta^*;\xi\cdot\eta)F(\eta)d\omega(\eta)- \Lambda^*_\xi\int_\Omega G(\Delta^*;\xi\cdot\eta)F(\eta)d\omega(\eta)\right|&=&0,
\end{eqnarray*}
where $\Lambda^*$ denotes one of the operators $\nabla^*$ or $L^*$.
\end{lem}

Relations for higher order derivatives are simple consequences of the above lemmas by use of well-known surface versions of Green's formulas that shift the differentiation from the convolution kernel to the convolved function $F$. Thus, they also require a higher smoothness of $F$.

\begin{cor}\label{cor:convbeltramig}
Let $G^\rho(\Delta^*; \cdot)$ be of class $ C^{(2)}([-1,1])$ and $F$ of class $C^{(1)}(\Omega)$. Then 
\begin{eqnarray*}
\lim_{\rho\to0+}\sup_{\xi\in\Omega}\left|\int_\Omega \Delta^*_\xi G^\rho(\Delta^*;\xi\cdot\eta)F(\eta)d\omega(\eta)- \Delta^*_\xi\int_\Omega G(\Delta^*;\xi\cdot\eta)F(\eta)d\omega(\eta)\right|&=&0.
\end{eqnarray*}
\end{cor}

\begin{cor}\label{cor:convtensbeltramig}
Let $G^\rho(\Delta^*; \cdot)$ be of class $ C^{(2)}([-1,1])$ and $f$ of class $c^{(1)}(\Omega)$. Then
\begin{eqnarray*}
\lim_{\rho\to0+}\sup_{\xi\in\Omega}&\!\!\!\!\bigg|&\!\!\!\!\int_\Omega \big(\big(\Lambda_1^{*}\big)_\xi\otimes \big(\Lambda_2^{*}\big)_\eta G^\rho(\Delta^*;\xi\cdot\eta)\big)f(\eta)d\omega(\eta)
\\&&-\big(\Lambda_1^{*}\big)_\xi\int_\Omega \big(\Lambda_2^{*}\big)_\eta G(\Delta^*;\xi\cdot\eta)\cdot f(\eta)d\omega(\eta)\bigg|=0,
\end{eqnarray*}
where $\Lambda_1^*$ and $\Lambda_2^*$ denote one of the operators $\nabla^*$ or $L^*$ ($\otimes$ denotes the tensor product $x\otimes y=xy^T$, for $x,y\in\mathbb{R}^3$).
\end{cor}

An adequate choice of $R^\rho$ admits an explicit statement on the convergence rate. More precisely, if $\int_{1-\rho}^1\left|R^\rho(t)\right|\,dt=\mathcal{O}(\rho)$ and $\int_{1-\rho}^1\left|\frac{d}{dt}R^\rho(t)\right|\,dt=\mathcal{O}(1)$, one can find
\begin{eqnarray*}
\left|\int_\Omega G^\rho(\Delta^*;\xi\cdot\eta)F(\eta)d\omega(\eta)-\int_\Omega G(\Delta^*;\xi\cdot\eta)F(\eta)d\omega(\eta)\right|&=&\mathcal{O}(\rho\ln(\rho)),
\\\left|\int_\Omega \Lambda^*_\xi G^\rho(\Delta^*;\xi\cdot\eta)F(\eta)d\omega(\eta)- \Lambda^*_\xi\int_\Omega G(\Delta^*;\xi\cdot\eta)F(\eta)d\omega(\eta)\right|&=&\mathcal{O}(\rho^{\frac{1}{2}}),
\end{eqnarray*}
for $F$ of class $C^{(0)}(\Omega)$. If $F$ is of class $C^{(1)}(\Omega)$, it even holds 
\begin{eqnarray*}
\left|\int_\Omega \Lambda^*_\xi G^\rho(\Delta^*;\xi\cdot\eta)F(\eta)d\omega(\eta)- \Lambda^*_\xi\int_\Omega G(\Delta^*;\xi\cdot\eta)F(\eta)d\omega(\eta)\right|&=&\mathcal{O}(\rho\ln(\rho)).
\end{eqnarray*}
The conditions on the regularization function are satisfied, e.g., by the choice of $R^\rho$ as the truncated Taylor series of $G(\Delta^*;\cdot)$.

\subsection{Single Layer Kernel}\label{sec:singlelayer}

By the singel layer kernel we denote the convolution kernel of the integral operator $D^{-1}$, with $D$ formally given as in (\ref{eqn:d}). Observing that $\Delta^*Y_{n,k}=-n(n+1)Y_{n,k}$, the fractional pseudodifferential operator $D$, mapping the Sobolev space $H_s(\Omega)$ into $H_{s-1}(\Omega)$, can be defined via
\begin{eqnarray}
DF=\left(-\Delta^*+\textnormal{\footnotesize $\frac{1}{4}$}\right)^{\frac{1}{2}}F=\sum_{n=0}^\infty\sum_{k=1}^{2n+1}\big(n+\textnormal{\footnotesize$\frac{1}{2}$}\big)\,(F,Y_{n,k})_{L^2(\Omega)}Y_{n,k},
\end{eqnarray}
for $F$ of class $H_s(\Omega)$, where $(\cdot,\cdot)_{L^2(\Omega)}$ denotes the inner product $(F,G)_{L^2(\Omega)}=\int_\Omega F(\eta)G(\eta)d\omega(\eta)$. Its inverse $D^{-1}$, mapping $H_{s-1}(\Omega)$ into $H_{s}(\Omega)$, is correspondingly given by 
\begin{eqnarray}
D^{-1}F=\left(-\Delta^*+\textnormal{\footnotesize $\frac{1}{4}$}\right)^{-\frac{1}{2}}F=\sum_{n=0}^\infty\sum_{k=1}^{2n+1}\frac{1}{n+\textnormal{\footnotesize$\frac{1}{2}$}}\,(F,Y_{n,k})_{L^2(\Omega)}Y_{n,k}, \label{eqn:singlelayerpseudodiff}
\end{eqnarray}
for $F$ of class $H_{s-1}(\Omega)$. From the addition theorem and the power series in Lemma \ref{lem:seriespotfunc}, it is easy to derive the integral representation
\begin{eqnarray*}
\sum_{n=0}^\infty\sum_{k=1}^{2n+1}\frac{1}{n+\frac{1}{2}}(F,Y_{n,k})_{L^2(\Omega)}Y_{n,k}(\xi)=\frac{1}{2\sqrt{2}\pi}\int_\Omega \frac{1}{\sqrt{1-\xi\cdot\eta}}F(\eta)d\omega(\eta),\quad \xi\in\Omega.
\end{eqnarray*}
The function
\begin{eqnarray}
S(\xi\cdot\eta)=\frac{1}{\sqrt{2}}\frac{1}{\sqrt{1-\xi\cdot\eta}},\quad 1-\xi\cdot\eta>0,
\end{eqnarray}
is from now on called the \emph{single layer kernel}, and denotes the starting point for our further considerations. The integral operator $D^{-1}$ is called the single layer operator. For a more general and detailed overview on spherical pseudodifferential operators and the definition of Sobolev spaces, the reader is referred to, e.g., \cite{free98} and \cite{svensson}. Since we are dealing with continuously differentiable functions in the remainder of this paper, it should be remarked that $D^{-1}$ actually maps $C^{(k)}(\Omega)$ into $C^{(k)}(\Omega)$, $k\in\mathbb{N}_0$.

In analogy to Green's function, one can define a spatial regularization of the single layer kernel.

\begin{defi}\label{def:regsingkernel}
Let $\rho>0$ and $R^\rho$ a non-negative function of class $C^{(n)}([-1,1])$, $n\in\mathbb{N}$ fixed, satisfying
\begin{eqnarray*}\label{eqn:propregsing}
\lim_{\rho\to0+}\rho^k\int_{1-\rho}^1\left|\left(\frac{d}{dt}\right)^kR^\rho(t)\right|\,dt=0,\quad k=0,1,
\end{eqnarray*}
and 
\begin{equation*}
\left[\left(\frac{d}{dt}\right)^k R^\rho(t)\right]_{t=1-\rho}=\left[\left(\frac{d}{dt}\right)^kS(t)\right]_{t=1-\rho},\quad k=0,1,\ldots,n. 
\end{equation*}
Then the function
\begin{eqnarray*}
S^\rho(\xi\cdot\eta)=\left\{\begin{array}{ll}
\frac{1}{\sqrt{2}}\frac{1}{\sqrt{1-\xi\cdot\eta}},&1-\xi\cdot\eta\geq\rho,
\\[1.25ex]R^\rho(\xi\cdot\eta),&1-\xi\cdot\eta<\rho,
\end{array}\right.
\end{eqnarray*}
is called \emph{regularized single layer kernel (of order $n$)}.
\end{defi}

The regularizing function $R^\rho$ is generally chosen as the Taylor series of $S$ centered at $1-\rho$ and truncated at degree $n$. An exemplary plot for different scaling parameters $\rho$ can be found in Figure \ref{fig:reggreenfunc}. The special cases of a linear or quadratic regularization have been applied to multiscale methods in physical geodesy, e.g., in \cite{freewolf09}, \cite{free09}. In the Euclidean space $\mathbb{R}^3$, the kernel $S$ can be related to the fundamental solution of the Laplace operator $\Delta$. A different kind of regularization for that kernel is treated in \cite{akram10}. In our setting, we obtain the following limit relation in the same manner as for the Green function case in the previous subsection.

\begin{lem}\label{lem:convs}
Let  $S^\rho$ be of class $ C^{(1)}([-1, 1])$ and $F$ of class $C^{(0)}(\Omega)$. Then we have
\begin{eqnarray*}
\lim_{\rho\to0+}\sup_{\xi\in\Omega}\left|\int_\Omega S^\rho(\xi\cdot\eta)F(\eta)d\omega(\eta)-\int_\Omega S(\xi\cdot\eta)F(\eta)d\omega(\eta)\right|=0.
\end{eqnarray*}
\end{lem}

For the relations involving the surface gradient and the surface curl gradient, one has to observe that $\eta\mapsto \nabla_\xi^*S(\xi\cdot\eta)$ and $\eta\mapsto L_\xi^*S(\xi\cdot\eta)$ are not integrable on the sphere $\Omega$. However, if $F$ is of class $C^{(1)}(\Omega)$ and $\mathbf{t}_\xi\in\mathbb{R}^{3\times3}$ denotes the rotation matrix with $\mathbf{t}_\xi\xi=\eps^3=(0,0,1)^T$, we obtain
\begin{eqnarray*}
\nabla_\xi^*\int_\Omega S(\xi\cdot\eta) F(\eta)d\omega(\eta)&=&\nabla_\xi^*\int_\Omega S(\xi\cdot\mathbf{t}_\xi^T\eta)F(\mathbf{t}_\xi^T\eta)d\omega(\eta)\nonumber
\\&=&\nabla_\xi^*\int_\Omega S(\eta_3)F(\mathbf{t}_\xi^T\eta)d\omega(\eta)
\\&=&\int_\Omega S(\eta_3)\nabla_\xi^*F(\mathbf{t}_\xi^T\eta)d\omega(\eta),\quad\xi\in\Omega,\nonumber
\end{eqnarray*}
where $\eta=(\eta_1,\eta_2,\eta_3)^T\in\Omega$. Furthermore, regularizing the single layer kernel yields
\begin{eqnarray*}
\int_\Omega \nabla_\xi^*S^\rho(\xi\cdot\eta) F(\eta)d\omega(\eta)&=&\nabla_\xi^*\int_\Omega S^\rho(\xi\cdot\eta) F(\eta)d\omega(\eta)
\\&=&\nabla_\xi^*\int_\Omega S^\rho(\eta_3)F(\mathbf{t}_\xi^T\eta)d\omega(\eta)
\\&=&\int_\Omega S^\rho(\eta_3)\nabla_\xi^*F(\mathbf{t}_\xi^T\eta)d\omega(\eta),\quad\xi\in\Omega,
\end{eqnarray*}
so that the previous lemma implies the desired relations. 

\begin{lem}\label{lem:convsgradsing}
Let $F$ be of class $C^{(1)}(\Omega)$ and $S^\rho$ of class $C^{(1)}([-1,1])$. Then we have
\begin{eqnarray*}
\lim_{\rho\to0+} \sup_{\xi\in\Omega} \left|\int_\Omega \Lambda^*_\xi S^\rho(\xi\cdot\eta) F(\eta)d\omega(\eta)-\Lambda_\xi^*\int_\Omega S(\xi\cdot\eta) F(\eta)d\omega(\eta)\right|&=&0,
\end{eqnarray*}
where $\Lambda^*$ denotes one of the operators $\nabla^*$ or $L^*$.
\end{lem}

Relations for higher order differential operators follow analogously. Of more interest to us are combinations of the single layer operator with Green's function for the Beltrami operator. It holds, e.g., that
\begin{eqnarray}\label{eqn:singgreenrel}
\lim_{\rho\to0+}\sup_{\xi\in\Omega}\left|\int_\Omega  \left(\Lambda_\xi^* D_\xi^{-1}G^\rho(\Delta^* ;\xi\cdot\eta)\right)\,F(\eta)d\omega(\eta)-\Lambda_\xi^* D_\xi^{-1}\int_\Omega G(\Delta^*;\xi\cdot\eta)F(\eta)d\omega(\eta)\right|=0.\,\,\,\,
\end{eqnarray}
However, it is difficult to explicitly calculate $D_\xi^{-1}G^\rho(\Delta^* ;\xi\cdot\eta)$, as it would be required for our later applications. $D_\xi^{-1}G(\Delta^* ;\xi\cdot\eta)$, on the other hand, can be calculated, and a regularization afterwards yields a similar limit relation.

\begin{lem}
For $\xi,\eta\in\Omega$ we have,
\begin{eqnarray*}
D_\xi^{-1}G(\Delta^*;\xi\cdot\eta)=\frac{1}{2\pi}\ln\left((1+\xi\cdot\eta)\left(\frac{1}{2}-\frac{1}{1-2S(\xi\cdot\eta)}\right)\right)-\frac{1}{2\pi}.
\end{eqnarray*}
\end{lem}

\begin{prf}
Lemma \ref{lem:genseriesgreenfunc} and the pseudodifferential representation (\ref{eqn:singlelayerpseudodiff}) imply
\begin{eqnarray*}
D_\xi^{-1}G(\Delta^*;\xi\cdot\eta)&=&\sum_{n=1}^\infty\frac{1}{n+\frac{1}{2}}\frac{2n+1}{4\pi}\frac{1}{-n(n+1)}P_n(\xi\cdot\eta)
\\&=&\frac{1}{2\pi}\sum_{n=1}^\infty\frac{1}{n+1}P_n(\xi\cdot\eta)-\frac{1}{2\pi}\sum_{n=1}^\infty\frac{1}{n}P_n(\xi\cdot\eta)
\\&=&\frac{1}{2\pi}\ln\left(1+\frac{\sqrt{2}}{\sqrt{1-\xi\cdot\eta}}\right)-\frac{1}{2\pi}\ln\left(\frac{\sqrt{2}\sqrt{1-\xi\cdot\eta}-1+\xi\cdot\eta}{1-(\xi\cdot\eta)^2}\right)
-\frac{1}{2\pi}(1+\ln\left(2\right))
\\&=&\frac{1}{2\pi}\ln\left((1+\xi\cdot\eta)\left(\frac{1}{2}-\frac{1}{1-2S(\xi\cdot\eta)}\right)\right)-\frac{1}{2\pi},\pagebreak[0]
\end{eqnarray*}\pagebreak[0]which is well-defined for every $\xi,\eta\in\Omega$.\pagebreak[4]
\end{prf}

The above derived representation implies that $(\xi,\eta)\mapsto D_\xi^{-1}G(\Delta^*;\xi\cdot\eta)$ is zonal and of class $C^{(1)}(\Omega\times\Omega)$. Some lengthy but basic computations yield
\begin{eqnarray*}
\nabla_\xi^*D_\xi^{-1}G(\Delta^*;\xi\cdot\eta)=\frac{1}{2\pi}\left(\frac{1}{2}-S(\xi\cdot\eta)-\frac{1}{2+4S(\xi\cdot\eta)}\right)(\eta-(\xi\cdot\eta)\xi),\quad \xi,\eta\in\Omega.
\end{eqnarray*}
A further application of the surface gradient causes a singularity of type $\mathcal{O}((1-\xi\cdot\eta)^{-\frac{1}{2}})$. Therefore, we do the following regularization for $\rho>0$,
\begin{eqnarray}\label{eqn:regsinggreen}
s_{\nabla^*}^\rho(\xi,\eta)=\frac{1}{2\pi}\left(\frac{1}{2}-S^\rho(\xi\cdot\eta)-\frac{1}{2+4S^\rho(\xi\cdot\eta)}\right)(\eta-(\xi\cdot\eta)\xi),\quad \xi,\eta\in\Omega.
\end{eqnarray}
For this kernel we can calculate
\begin{eqnarray*}
\nabla_\xi^*\otimes s_{\nabla^*}^\rho(\eta,\xi)&=&\frac{1}{2\pi}\left(\frac{1}{2}-S^\rho(\xi\cdot\eta)-\frac{1}{2+4S^\rho(\xi\cdot\eta)}\right)\nabla_\xi^*\otimes(\xi-(\xi\cdot\eta)\eta)
\\&&+\frac{1}{2\pi}\left(-\big(S^\rho\big)'(\xi\cdot\eta)+\frac{4\big(S^\rho\big)'(\xi\cdot\eta)}{(2+4S^\rho(\xi\cdot\eta))^2}\right)(\eta-(\xi\cdot\eta)\xi)\otimes (\xi-(\xi\cdot\eta)\eta),
\end{eqnarray*}
where $\big(S^\rho\big)'$ denotes the one-dimensional derivative of $S^\rho$. The analogous procedure works for the surface curl gradient, and we have for $\rho>0$ that
\begin{eqnarray}\label{eqn:regsinggreen2}
s_{L^*}^\rho(\xi,\eta)=\frac{1}{2\pi}\left(\frac{1}{2}-S^\rho(\xi\cdot\eta)-\frac{1}{2+4S^\rho(\xi\cdot\eta)}\right)(\xi\wedge\eta), \quad \xi,\eta\in\Omega,
\end{eqnarray}
and 
\begin{eqnarray*}
L_\xi^*\otimes s_{L^*}^\rho(\eta,\xi)&=&\frac{1}{2\pi}\left(\frac{1}{2}-S^\rho(\xi\cdot\eta)-\frac{1}{2+4S^\rho(\xi\cdot\eta)}\right)L_\xi^*\otimes(\eta\wedge\xi)
\\&&+\frac{1}{2\pi}\left(-\big(S^\rho\big)'(\xi\cdot\eta)+\frac{4\big(S^\rho\big)'(\xi\cdot\eta)}{(2+4S^\rho(\xi\cdot\eta))^2}\right)(\xi\wedge\eta)\otimes (\eta\wedge\xi).
\end{eqnarray*}

Thus, relation (\ref{eqn:singgreenrel}) can be formulated in the following numerically more advantageous way.

\begin{lem}\label{lem:limsinggreen}
Let $F$ be of class $C^{(0)}(\Omega)$ and $S^\rho$ of class $C^{(1)}([-1,1])$. Then we get with $s_{\nabla^*}^\rho(\cdot,\cdot)$ and $s_{L^*}^\rho(\cdot,\cdot)$ as in (\ref{eqn:regsinggreen}) and (\ref{eqn:regsinggreen2}), respectively, that
\begin{eqnarray*}
\lim_{\rho\to0+}\sup_{\xi\in\Omega}\left|\int_\Omega s_{\Lambda^*}^\rho(\xi,\eta)F(\eta)\,d\omega(\eta)-\Lambda_\xi^*D_\xi^{-1}\int_\Omega G(\Delta^*;\xi\cdot\eta) F(\eta)\,d\omega(\eta)\right|&=&0,
\end{eqnarray*}
where $\Lambda^*$ denotes one of the operators $\nabla^*$ or $L^*$.
\end{lem}

The relation we are actually aiming at, and which we require in later applications, is the following tensorial one.

\begin{lem}\label{lem:limdiffsinggreen }
Let $f$ be of class $c^{(1)}(\Omega)$ and $S^\rho$ of class $C^{(1)}([-1,1])$. Then we get with $s_{\nabla^*}^\rho(\cdot,\cdot)$ and $s_{L^*}^\rho(\cdot,\cdot)$ as in (\ref{eqn:regsinggreen}) and (\ref{eqn:regsinggreen2}), respectively, that
\begin{eqnarray*}
\lim_{\rho\to0+}\sup_{\xi\in\Omega}\left|\int_\Omega \left(\Lambda_\xi^*\otimes s_{\Lambda^*}^\rho(\eta,\xi)\right) f(\eta)\,d\omega(\eta)-\Lambda_\xi^*\int_\Omega \left(\Lambda_\eta^*D_\xi^{-1}G(\Delta^*;\xi\cdot\eta)\right)\cdot f(\eta)\,d\omega(\eta)\right|&=&0,
\end{eqnarray*}
where $\Lambda^*$ is one of the operators $\nabla^*$ or $L^*$.
\end{lem}

\begin{prf}
Since $|\nabla_\xi^*\otimes (\xi-(\xi\cdot\eta)\eta)|$ and $|f(\eta)|$ are uniformly bounded with respect to $\xi,\eta\in\Omega$ by some constant $M>0$, we get the following estimate for $\xi\in\Omega$ and $\rho>0$,
\begin{eqnarray}
&&\!\!\!\!\!\!\!\!\!\!\left|\int_\Omega \left(\nabla_\xi^*\otimes s_{\nabla^*}^\rho(\eta,\xi)\right) f(\eta)\,d\omega(\eta)-\int_\Omega \left(\nabla_\xi^*\otimes\nabla_\eta^*D_\xi^{-1}G(\Delta^*;\xi\cdot\eta)\right) f(\eta)\,d\omega(\eta)\right|\label{eqn:tenslimgreensing}
\\&\leq&\int_{\eta\in\Omega\atop 1-\xi\cdot\eta\leq \rho}\left|S(\xi\cdot\eta)+\frac{1}{2+4S(\xi\cdot\eta)}-S^\rho(\xi\cdot\eta)-\frac{1}{2+4S^\rho(\xi\cdot\eta)}\right|\nonumber
\\&&\qquad\qquad\times\left|\nabla_\xi^*\otimes (\xi-(\xi\cdot\eta)\eta)\right|\,|f(\eta)|\,d\omega(\eta)\nonumber
\\&&+\int_{\eta\in\Omega\atop 1-\xi\cdot\eta\leq \rho}\left|S(\xi\cdot\eta)^3-\frac{4S(\xi\cdot\eta)^3}{(2+4S(\xi\cdot\eta\nonumber))^2} -\big(S^\rho\big)'(\xi\cdot\eta)+\frac{4\big(S^\rho\big)'(\xi\cdot\eta)}{(2+4S^\rho(\xi\cdot\eta))^2}\right|
\\&&\qquad\qquad\quad\times\left|(\eta-(\xi\cdot\eta)\xi)\otimes (\xi-(\xi\cdot\eta)\eta)\right|\,|f(\eta)|\,d\omega(\eta)\nonumber
\\[1.25ex]&\leq& M^2 \int_{\eta\in\Omega\atop 1-\xi\cdot\eta\leq \rho}\left|S(\xi\cdot\eta)+\frac{1}{2+4S(\xi\cdot\eta)}-S^\rho(\xi\cdot\eta)-\frac{1}{2+4S^\rho(\xi\cdot\eta)}\right|\,d\omega(\eta)\nonumber
\\&&+\,M\int_{\eta\in\Omega\atop 1-\xi\cdot\eta\leq \rho}\left|S(\xi\cdot\eta)^3-\frac{4S(\xi\cdot\eta)^3}{(2+4S(\xi\cdot\eta))^2} -\big(S^\rho\big)'(\xi\cdot\eta)+\frac{4\big(S^\rho\big)'(\xi\cdot\eta)}{(2+4S^\rho(\xi\cdot\eta))^2}\right|\nonumber
\\&&\qquad\qquad\qquad\times|\eta-(\xi\cdot\eta)\xi|\,|\xi-(\xi\cdot\eta)\eta|d\omega(\eta).\nonumber
\end{eqnarray}
Observing 
\[
|\eta-(\xi\cdot\eta)\xi|\,|\xi-(\xi\cdot\eta)\eta|=\frac{1}{2S(\xi\cdot\eta)^2}(1+\xi\cdot\eta),
\]
the integrability of $\eta\mapsto S(\xi\cdot\eta)$ on the sphere $\Omega$, and the properties for $S^\rho$ from Definition \ref{def:regsingkernel}, we see that the integrals above vanish as $\rho$ tends to zero. Due to the zonality of the kernels, this convergence is uniform with respect to $\xi\in\Omega$. Furthermore, the convergence of (\ref{eqn:tenslimgreensing}) to zero additionally yields
\begin{eqnarray*}
\int_\Omega \left(\nabla_\xi^*\otimes\nabla_\eta^*D_\xi^{-1}G(\Delta^*;\xi\cdot\eta)\right) f(\eta)\,d\omega(\eta)=\nabla_\xi^*\int_\Omega \left(\nabla_\eta^*D_\xi^{-1}G(\Delta^*;\xi\cdot\eta)\right)\cdot f(\eta)\,d\omega(\eta),
\end{eqnarray*}
and therefore, the desired statement. The assertion for the surface curl gradient follows analogously.
\end{prf}

In analogy to the Green's function case, an adequate choice of $R^\rho$ admits an explicit statement on the convergence rate. More precisely, if $\int_{1-\rho}^1\left|R^\rho(t)\right|\,dt=\mathcal{O}(\rho)$ and $\int_{1-\rho}^1\left|\frac{d}{dt}R^\rho(t)\right|\,dt=\mathcal{O}(1)$, one obtains
\begin{eqnarray*}
\left|\int_\Omega \Lambda^*_\xi S^\rho(\xi\cdot\eta) F(\eta)d\omega(\eta)-\Lambda_\xi^*\int_\Omega S(\xi\cdot\eta) F(\eta)d\omega(\eta)\right|&=&\mathcal{O}(\rho^{\frac{1}{2}}).
\end{eqnarray*}
For the combination of Green's function and the single layer kernel the same convergence rate holds true,
\begin{eqnarray*}
\left|\int_\Omega \left(\Lambda_\xi^*\otimes s_{\Lambda^*}^\rho(\eta,\xi)\right) f(\eta)\,d\omega(\eta)-\Lambda_\xi^*\int_\Omega \left(\Lambda_\eta^*D_\xi^{-1}G(\Delta^*;\xi\cdot\eta)\right)\cdot f(\eta)\,d\omega(\eta)\right|&=&\mathcal{O}(\rho^{\frac{1}{2}}),
\end{eqnarray*}
for $F$ of class $C^{(1)}(\Omega)$ and $f$ of class $c^{(1)}(\Omega)$. The conditions on the regularization function are satisfied, e.g., by the choice of $R^\rho$ as the truncated Taylor series of $S$.

\section{Spherical Decompositions}\label{sec:decomps}

We introduce two decompositions relating to the two sets of vector spherical harmonics from Section \ref{sec:basics}. The first one, relating to $y^{(i)}_{n,k}$ and the operators $o^{(i)}$, respectively, is the well-known spherical Helmholtz decomposition. It decomposes a vector field into its radial part and two tangential parts. The second one, relating to $\tilde{y}_{n,k}^{(i)}$ and the operators $\tilde{o}^{(i)}$, respectively, is crucial for the separation with respect to the sources and is presented in some detail in this section. A more general overview on similar spherical decompositions can be found, e.g., in \cite{gerhards10}. 
 
\begin{thm}[Helmholtz decomposition]\label{thm:helmholtzglob}
Let $f$ be of class $c^{(1)}(\Omega)$. Then there exist uniquely defined scalar fields $F_1$ of class $C^{(1)}(\Omega)$ and $F_2,F_3$ of class $C^{(2)}(\Omega)$ satisfying
\[
\frac{1}{4\pi}\int_{\Omega}F_i(\eta)d\omega(\eta)=0,\quad i=2,3,
\]
such that
\begin{eqnarray*}
f(\xi)&=&o^{(1)}_\xi F_1(\xi)+o^{(2)}_\xi F_2(\xi)+o^{(3)}_\xi F_3(\xi),\quad \xi\in\Omega.
\end{eqnarray*}
\end{thm}

A proof of this decomposition can be found, e.g., in \cite{back96}. Using Green's function for the Beltrami operator and Theorem \ref{thm:fundthmgradcurl} yields the following representations for the Helmholtz scalars,
\begin{eqnarray}
F_1(\xi)&=&\xi\cdot f(\xi),\quad\xi\in\Omega,\label{eqn:h1}
\\F_2(\xi)&=&-\int_\Omega\big(\nabla_\eta^*G(\Delta^*;\xi\cdot\eta)\big)\cdot f(\eta)\,d\omega(\eta),\quad\xi\in\Omega,\label{eqn:hemlf2rep}
\\F_3(\xi)&=&-\int_\Omega\big(L_\eta^*G(\Delta^*;\xi\cdot\eta)\big)\cdot f(\eta)\,d\omega(\eta),\quad\xi\in\Omega.\label{eqn:hemlf3rep}
\end{eqnarray}

If $F_1$ additionally has vanishing integral mean value, i.e., $\int_\Omega F_1(\eta)d\omega(\eta)=0$ (as is the case for functions satisfying the pre-Maxwell equations), there exists a function $U$ of class $C^{(2)}(\Omega)$ with $\Delta^*U=F_1$, such that Theorem \ref{thm:fundthmdelta} implies
\begin{eqnarray}
F_1(\xi)&=&\Delta_\xi^*\int_\Omega G(\Delta^*;\xi\cdot\eta)\eta\cdot f(\eta)d\omega(\eta),\quad\xi\in\Omega.\label{eqn:h4}
\end{eqnarray}

While the orthogonality is the main property of the Helmholtz decomposition, a representation with respect to the operators $\tilde{o}^{(i)}$ is of special interest in geomagnetic modeling. In order to obtain a representation of the corresponding scalars, we rewrite (\ref{eqn:to1operator})--(\ref{eqn:to3operator}) as
\begin{eqnarray}
o^{(1)}&=&\textnormal{\footnotesize $\frac{1}{2}$}\tilde{o}^{(1)} D^{-1}+\textnormal{\footnotesize $\frac{1}{2}$}\tilde{o}^{(2)} D^{-1},\label{eqn:tto1operator}
\\o^{(2)}&=&\textnormal{\footnotesize $\frac{1}{2}$}\tilde{o}^{(1)} \left(\textnormal{\footnotesize $\frac{1}{2}$}D^{-1}-1\right)+\textnormal{\footnotesize $\frac{1}{2}$}\tilde{o}^{(2)} \left(\textnormal{\footnotesize $\frac{1}{2}$}D^{-1}+1\right),\label{eqn:tto2operator}
\\o^{(3)}&=&\tilde{o}^{(3)}.\label{eqn:tto3operator}
\end{eqnarray}
This gives us the necessary representation to prove the following decomposition theorem.

\begin{thm}\label{thm:althelmholtz}
Let $f$ be of class $c^{(1)}(\Omega)$. Then there exist uniquely defined scalar fields $\tilde{F}_1,\tilde{F}_2$ of class $C^{(1)}(\Omega)$ and $\tilde{F}_3$ of class $C^{(2)}(\Omega)$ satisfying
\begin{eqnarray*}
\frac{1}{4\pi}\int_{\Omega}\tilde{F}_3(\eta)d\omega(\eta)&=&0,
\\\frac{1}{4\pi}\int_{\Omega}\tilde{F}_1(\eta)-\tilde{F}_2(\eta)d\omega(\eta)&=&0,
\end{eqnarray*}
such that
\begin{eqnarray*}
f(\xi)&=&\tilde{o}^{(1)}_\xi \tilde{F}_1(\xi)+\tilde{o}^{(2)}_\xi \tilde{F}_2(\xi)+\tilde{o}^{(3)}_\xi \tilde{F}_3(\xi),\quad \xi\in\Omega.
\end{eqnarray*}
The scalars $\tilde{F}_1,\tilde{F}_2$ and $\tilde{F}_3$ can be represented by
\begin{eqnarray}
\tilde{F}_1&=&\textnormal{\footnotesize $\frac{1}{2}$}D^{-1}F_1+\textnormal{\footnotesize $\frac{1}{4}$}D^{-1}F_2-\textnormal{\footnotesize $\frac{1}{2}$}F_2,\label{eqn:tf1}
\\\tilde{F}_2&=&\textnormal{\footnotesize $\frac{1}{2}$}D^{-1}F_1+\textnormal{\footnotesize $\frac{1}{4}$}D^{-1}F_2+\textnormal{\footnotesize $\frac{1}{2}$}F_2,\label{eqn:tf2}
\\\tilde{F}_3&=&F_3,\label{eqn:tf3}
\end{eqnarray}
with $F_1,F_2$ and $F_3$ being the uniquely determined functions of the Helmholtz decomposition in Theorem \ref{thm:helmholtzglob}.
\end{thm}

\begin{prf}
Applying the Helmholtz decomposition to $f$ and using (\ref{eqn:tto1operator})--(\ref{eqn:tto3operator}), we get on $\Omega$,
\begin{eqnarray*}
f&=&o^{(1)}F_1+o^{(2)}F_2+o^{(3)}F_3
\\&=&\frac{1}{2}\tilde{o}^{(1)}D^{-1}F_1+\frac{1}{2}\tilde{o}^{(2)}D^{-1}F_1+\frac{1}{2}\tilde{o}^{(1)}\left(\frac{1}{2}D^{-1}-1\right)F_2+\frac{1}{2}\tilde{o}^{(2)}\left(\frac{1}{2}D^{-1}+1\right)F_2+\tilde{o}^{(3)}F_3
\\&=&\tilde{o}^{(1)}\left(\frac{1}{2}D^{-1}F_1+\frac{1}{4}D^{-1}F_2-\frac{1}{2}F_2\right)+\tilde{o}^{(2)}\left(\frac{1}{2}D^{-1}F_1+\frac{1}{4}D^{-1}F_2+\frac{1}{2}F_2\right)+\tilde{o}^{(3)}F_3.
\end{eqnarray*}
This implies a decomposition as stated in the theorem. Due to the uniqueness of the Helmholtz representation, it follows directly that $\tilde{F}_3$ is defined uniquely when having a vanishing integral mean value. For the uniqueness of $\tilde{F}_1$ and $\tilde{F}_2$ it is sufficient to show that $f(\xi)=0$, $\xi\in\Omega$, only has the trivial decomposition with respect to the operators $\tilde{o}^{(i)}$, $i=1,2,3$. If
\[
\tilde{o}^{(1)}_\xi\tilde{F}_1(\xi)+\tilde{o}^{(2)}_\xi\tilde{F}_2(\xi)+\tilde{o}^{(3)}_\xi\tilde{F}_3(\xi)=0, \quad\xi\in\Omega,
\]
we get from (\ref{eqn:to1operator})--(\ref{eqn:to3operator}) that
\[
o^{(1)}_\xi\left(\left(D_\xi+\textnormal{\footnotesize $\frac{1}{2}$}\right)\tilde{F}_1(\xi)+\left(D_\xi-\textnormal{\footnotesize $\frac{1}{2}$}\right)\tilde{F}_2(\xi)\right)+o^{(2)}_\xi\left(\tilde{F}_2(\xi)-\tilde{F}_1(\xi)\right)+o^{(3)}_\xi\tilde{F}_3(\xi)=0,\quad\xi\in\Omega.
\]
The uniqueness of the Helmholtz decomposition then implies
\begin{eqnarray*}
\tilde{F}_2(\xi)-\tilde{F}_1(\xi)&=&0,\quad\xi\in\Omega,
\\\left(D_\xi+\textnormal{\footnotesize $\frac{1}{2}$}\right)\tilde{F}_1(\xi)+\left(D_\xi-\textnormal{\footnotesize $\frac{1}{2}$}\right)\tilde{F}_2(\xi)&=&0,\quad \xi\in\Omega,
\end{eqnarray*}
if $\frac{1}{4\pi}\int_{\Omega}\tilde{F}_1(\eta)-\tilde{F}_2(\eta)d\omega(\eta)=0$, which gives us
\[
D_\xi\tilde{F}_1(\xi)=0,\quad\xi\in\Omega.
\]
Thus, $\tilde{F}_1(\xi)=0$, $\xi\in\Omega$, since $D$ is injective, and it follows $\tilde{F}_2(\xi)=0$, $\xi\in\Omega$, so that uniqueness is given for this  decomposition.
\end{prf}

The theorem above yields, by use of (\ref{eqn:h1})--(\ref{eqn:hemlf3rep}), a representation of the scalars $\tilde{F}_i$. Of importance in the later application, however, are the vectorial quantities $\tilde{o}^{(i)}\tilde{F}_i$. Thus, we first calculate from (\ref{eqn:tf1}) that
\begin{eqnarray}
\tilde{o}^{(1)}\tilde{F}_1(\xi)&=&\frac{1}{2}\xi(\xi\cdot f(\xi))+\frac{1}{4}\xi D_\xi^{-1}(\xi\cdot f(\xi))
-\frac{1}{8}\xi D_\xi^{-1}\int_\Omega\nabla_\eta^*G(\Delta^*;\xi\cdot\eta)\cdot f(\eta)d\omega(\eta)\nonumber
\\&&+\frac{1}{2}\xi D_\xi\int_\Omega\nabla_\eta^*G(\Delta^*;\xi\cdot\eta)\cdot f(\eta)d\omega(\eta)\label{eqn:tf1d}
\\&&-\frac{1}{2}\nabla_\xi^*D_\xi^{-1}(\xi\cdot f(\xi)) +\frac{1}{4}\nabla_\xi^*D_\xi^{-1}\int_\Omega\nabla_\eta^*G(\Delta^*;\xi\cdot\eta)\cdot f(\eta)d\omega(\eta)\nonumber
\\&&-\frac{1}{2}\nabla_\xi^*\int_\Omega\nabla_\eta^*G(\Delta^*;\xi\cdot\eta)\cdot f(\eta)d\omega(\eta),\qquad\xi\in\Omega.\nonumber
\end{eqnarray}
The expression in the second row, involving the operator $D$, is unfortunate since we have no explicit representation for the corresponding regularized convolution kernel. Observing 
\[
D=D^{-1}\Big(-\Delta^*+\frac{1}{4}\Big),
\]
this can be circumvented by rewriting
\begin{eqnarray*}
&&\!\!\!\!\!\!\!\!\!\!\!\!\!\!D_\xi\int_\Omega \nabla_\eta^*G(\Delta^*;\xi\cdot\eta) \cdot f(\eta)d\omega(\eta)
\\&=&\frac{1}{4}D_\xi^{-1}\int_\Omega\nabla_\eta^* G(\Delta^*;\xi\cdot\eta) \cdot f(\eta)d\omega(\eta)+D_\xi^{-1}\Delta_\xi^*\int_\Omega G(\Delta^*;\xi\cdot\eta) \nabla_\eta^*\cdot f_{tan}(\eta)d\omega(\eta)
\\&=&\frac{1}{4}D_\xi^{-1}\int_\Omega\nabla_\eta^* G(\Delta^*;\xi\cdot\eta) \cdot f(\eta)d\omega(\eta)+D_\xi^{-1}\nabla_\xi^*\cdot f_{tan}(\xi),\qquad\xi\in\Omega,
\end{eqnarray*}
where Theorem \ref{thm:fundthmdelta} has been used in the last step, assuming $f$ to be of class $c^{(2)}(\Omega)$. Application of the above to (\ref{eqn:tf1d}) provides an easier calculable expression
\begin{eqnarray*}
\tilde{o}^{(1)}\tilde{F}_1(\xi)&=&\frac{1}{2}\xi(\xi\cdot f(\xi))+\frac{1}{4}\xi D_\xi^{-1}(\xi\cdot f(\xi))+\frac{1}{2}\xi D_\xi^{-1}\nabla_\xi^*\cdot f_{tan}(\xi)-\frac{1}{2}\nabla_\xi^*D_\xi^{-1}(\xi\cdot f(\xi))
\\&&+\frac{1}{4}\nabla_\xi^*\int_\Omega\nabla_\eta^*D_\xi^{-1}G(\Delta^*;\xi\cdot\eta)\cdot f(\eta)d\omega(\eta)-\frac{1}{2}\nabla_\xi^*\int_\Omega\nabla_\eta^*G(\Delta^*;\xi\cdot\eta)\cdot f(\eta)d\omega(\eta).
\end{eqnarray*}
Since the occurring differential operators and the integration cannot be interchanged without restriction, we need to switch to the regularized versions of the single layer kernel and the Green function for the Beltrami operator. Then it is valid to set
\begin{eqnarray}
\tilde{f}^{(1)}_\rho(\xi)&=&\frac{1}{2}\xi(\xi\cdot f(\xi))+\frac{1}{8\pi}\xi \int_\Omega S^\rho(\xi\cdot\eta) \,\eta\cdot f(\eta)d\omega(\eta)-\frac{1}{4\pi}\xi \int_\Omega \nabla_\eta^*S^\rho(\xi\cdot\eta)\cdot f(\eta)d\omega(\eta)\nonumber
\\&&-\frac{1}{4\pi}\int_\Omega \nabla_\xi^*S^\rho(\xi\cdot\eta)\,\eta\cdot f(\eta)d\omega(\eta)+\frac{1}{4}\int_\Omega\left(\nabla_\xi^*\otimes s_{\nabla^*}^\rho(\eta,\xi)\right) f(\eta)d\omega(\eta)\label{eqn:tf1approx}
\\&&-\frac{1}{2}\int_\Omega\left(\nabla_\xi^*\otimes\nabla_\eta^*G^\rho(\Delta^*;\xi\cdot\eta)\right) f(\eta)d\omega(\eta),\qquad\xi\in\Omega,\nonumber
\end{eqnarray}
for $\rho>0$. If $G^\rho(\Delta^*;\cdot)$ is of class $C^{(2)}([-1,1])$ and $S^\rho$ of class $C^{(1)}([-1,1])$, the considerations in Section \ref{sec:kernels} imply
\begin{eqnarray}
\lim_{\rho\to0+}\sup_{\xi\in\Omega}\left|\tilde{f}^{(1)}_\rho(\xi)-\tilde{o}^{(1)}\tilde{F}_1(\xi)\right|=0.\label{eqn:limittf1rho}
\end{eqnarray}

Analogous computations for $\tilde{o}^{(2)}\tilde{F}_2$ yield a regularization
\begin{eqnarray}
\tilde{f}^{(2)}_\rho(\xi)&=&\frac{1}{2}\xi(\xi\cdot f(\xi))-\frac{1}{8\pi}\xi \int_\Omega S^\rho(\xi\cdot\eta) \,\eta\cdot f(\eta)d\omega(\eta)+\frac{1}{4\pi}\xi \int_\Omega \nabla_\eta^*S^\rho(\xi\cdot\eta)\cdot f(\eta)d\omega(\eta)\nonumber
\\&&+\frac{1}{4\pi}\int_\Omega \nabla_\xi^*S^\rho(\xi\cdot\eta)\,\eta\cdot f(\eta)d\omega(\eta)-\frac{1}{4}\int_\Omega\left(\nabla_\xi^*\otimes s_{\nabla^*}^\rho(\eta,\xi)\right) f(\eta)d\omega(\eta)\label{eqn:tf2approx}
\\&&-\frac{1}{2}\int_\Omega\left(\nabla_\xi^*\otimes\nabla_\eta^*G^\rho(\Delta^*;\xi\cdot\eta)\right) f(\eta)d\omega(\eta),\qquad\xi\in\Omega,\nonumber
\end{eqnarray}
and
\begin{eqnarray}
\lim_{\rho\to0+}\sup_{\xi\in\Omega}\left|\tilde{f}^{(2)}_\rho(\xi)-\tilde{o}^{(2)}\tilde{F}_2(\xi)\right|=0.\label{eqn:limittf2rho}
\end{eqnarray}
Finally, the determination of $\tilde{o}^{(3)}\tilde{F}_3$ corresponds to the calculation of the toroidal part of $f$. From (\ref{eqn:tf3}) and (\ref{eqn:hemlf3rep}), we get
\[
\tilde{o}^{(3)}\tilde{F}_3(\xi)=-L_\xi^*\int_\Omega L_\eta^*G(\Delta^*;\xi\cdot\eta)\cdot f(\eta)d\omega(\eta),\quad\xi\in\Omega.
\]
Corollary \ref{cor:convtensbeltramig} then implies for the regularized version 
\begin{eqnarray}
\tilde{f}^{(3)}_\rho(\xi)=-\int_\Omega\left(L_\xi^*\otimes L_\eta^*G^\rho(\Delta^*;\xi\cdot\eta)\right)f(r\eta)\,d\omega(\eta),\quad\xi\in\Omega,\label{eqn:tf3approx}
\end{eqnarray}
that for $G^\rho(\Delta^*;\cdot)$ of class $C^{(2)}([-1,1])$,
\begin{equation}
\lim_{\rho\to0+}\sup_{\xi\in\Omega}\left|\tilde{f}^{(3)}_\rho(\xi)-\tilde{o}^{(3)}\tilde{F}_3(\xi)\right|=0.\label{eqn:limittf3rho}
\end{equation}
Summarizing, we can state the following theorem.

\begin{thm}\label{lem:althelmapprox}
Let $f$ be of class $c^{(2)}(\Omega)$, with $\int_\Omega\eta\cdot f(\eta)d\omega(\eta)=0$. Furthermore, let the regularized Green function $G^\rho(\Delta^*;\cdot)$ be of class $C^{(2)}([-1,1])$ and the single layer kernel $S^\rho$ of class $C^{(1)}([-1,1])$. Then
\begin{eqnarray*}
f(\xi)=\tilde{o}^{(1)}\tilde{F}_1(\xi)+\tilde{o}^{(2)}\tilde{F}_2(\xi)+\tilde{o}^{(3)}\tilde{F}_3(\xi),\quad\xi\in\Omega,
\end{eqnarray*}
with 
\begin{equation*}
\lim_{\rho\to0+}\sup_{\xi\in\Omega}\left|\int_\Omega\boldsymbol{\Phi}_\rho^{(i)}(\xi,\eta)f(\eta)d\omega(\eta)-\tilde{o}^{(i)}\tilde{F}_i(\xi)\right|=0,
\end{equation*}
for $i=1,2,3$. The convolution kernels are given by 
\begin{eqnarray*}
\boldsymbol{\Phi}_\rho^{(1)}(\xi,\eta)&=&\xi\otimes\eta\left(\frac{1}{2}\Delta^*_\xi G^{\rho}(\Delta^*;\xi\cdot\eta)+\frac{1}{8\pi}S^{\rho}(\xi\cdot\eta)\right)-\frac{1}{4\pi}\xi\otimes\nabla_\eta^*S^{\rho}(\xi\cdot\eta)
\\&&+\frac{1}{4}\nabla_\xi^*\otimes s_{\nabla^*}^{\rho}(\eta,\xi)-\frac{1}{4\pi}\nabla_\xi^*S^{\rho}(\xi,\eta)\otimes\eta-\frac{1}{2}\nabla_\xi^*\otimes\nabla_\eta^*G^{\rho}(\Delta^*;\xi\cdot\eta),\quad \xi,\eta\in\Omega,
\\\boldsymbol{\Phi}_\rho^{(2)}(\xi,\eta)&=&\xi\otimes\eta\left(\frac{1}{2}\Delta^*_\xi G^{\rho}(\Delta^*;\xi\cdot\eta)-\frac{1}{8\pi}S^{\rho}(\xi\cdot\eta)\right)+\frac{1}{4\pi}\xi\otimes\nabla_\eta^*S^{\rho}(\xi\cdot\eta)
\\&&-\frac{1}{4}\nabla_\xi^*\otimes s_{\nabla^*}^{\rho}(\eta,\xi)+\frac{1}{4\pi}\nabla_\xi^*S^{\rho}(\xi,\eta)\otimes\eta-\frac{1}{2}\nabla_\xi^*\otimes\nabla_\eta^*G^{\rho}(\Delta^*;\xi\cdot\eta),\quad \xi,\eta\in\Omega,
\\[1.25ex]\boldsymbol{\Phi}_\rho^{(3)}(\xi,\eta)&=&-L_\xi^*\otimes L_\eta^*G^{\rho}(\Delta^*;\xi\cdot\eta),\quad \xi,\eta\in\Omega.
\end{eqnarray*}
\end{thm}

\begin{prf}
Since $\int_\Omega\eta\cdot f(\eta)d\omega(\eta)=\int_\Omega F_1(\eta)d\omega(\eta)=0$, representation (\ref{eqn:h4}) implies 
\begin{eqnarray*}
\frac{1}{2}\xi(\xi\cdot f(\xi))&=&\frac{1}{2}\xi\Delta_\xi^*\int_\Omega G(\Delta^*;\xi\cdot\eta)\eta\cdot f(\eta),\quad\xi\in\Omega.
\end{eqnarray*}
Substituting the Green function by its regularized counterpart, we obtain
\begin{eqnarray*}
\frac{1}{2}\xi\Delta_\xi^*\int_\Omega G^\rho(\Delta^*;\xi\cdot\eta)\eta\cdot f(\eta)\,d\omega(\eta)&=&\frac{1}{2}\xi\int_\Omega \Delta_\xi^*G^\rho(\Delta^*;\xi\cdot\eta)\eta\cdot f(\eta)\,d\omega(\eta)
\\&=&\int_\Omega \left(\xi\otimes\eta\left(\frac{1}{2}\Delta_\xi^*G^\rho(\Delta^*;\xi\cdot\eta)\right)\right) f(\eta)\,d\omega(\eta),\quad\xi\in\Omega.
\end{eqnarray*}
Analogously, the remaining integral expressions in (\ref{eqn:tf1approx}) can be written in terms of convolutions with tensorial kernels if this is not already the case. This yields the kernel $\boldsymbol{\Phi}_\rho^{(1)}(\cdot,\cdot)$ for $\tilde{f}_\rho^{(1)}$. Corollary \ref{cor:convbeltramig} and (\ref{eqn:limittf1rho}) provide the desired limit relation for $\tilde{o}^{(1)}\tilde{F}_1$. The same holds true for $\tilde{o}^{(2)}\tilde{F}_2$ and $\tilde{o}^{(3)}\tilde{F}_3$.
\end{prf}

\section{Multiscale Representation for the Separation of Sources}\label{sec:multiscale}

The kernels from Theorem \ref{lem:althelmapprox} are the main ingredient to the upcoming multiscale representation. They actually denote the so-called scaling kernels, while the differences for different parameters $\rho$ denote the corresponding wavelet kernels. But before we go into detail, we briefly want to motivate why the decomposition with respect to the operators $\tilde{o}^{(i)}$ can be called a separation with respect to the sources.

From now on, we denote the magnetic field by $b$ of class $c^{(2)}(\mathbb{R}^3)$, the corresponding source current density by $j$ of class $c^{(1)}(\mathbb{R}^3)$, and by $\mu_0$ we mean the vacuum permeability. Furthermore, we assume the pre-Maxwell equations
\begin{eqnarray*}
\nabla_x\wedge b(x)&=&\mu_0j(x),\quad x\in\mathbb{R}^3,
\\\nabla_x\cdot b(x)&=&0,\quad x\in\mathbb{R}^3,
\end{eqnarray*}
to be satisfied. As mentioned in the introduction, the Mie decomposition (see, e.g., \cite{back86} and \cite{back96}) yields poloidal fields $p_b$, $p_j$ and toroidal fields $q_b$, $q_j$, such that
\begin{eqnarray*}
b(x)&=&p_b(x)+q_b(x),\quad x\in\mathbb{R}^3,
\\j(x)&=&p_j(x)+q_j(x),\quad x\in\mathbb{R}^3.
\end{eqnarray*}
Making use of the law of Biot-Savart (see, e.g., \cite{jackson75}) and the fact that the poloidal magnetic field $p_b$ is solely produced by tangential toroidal current densities $q_j$, the poloidal magnetic field can be split up as follows,
\begin{eqnarray*}
p_b(x)=p_b^{int}(R;x)+p_b^{ext}(R;x),\quad x\in\mathbb{R}^3\setminus\Omega_R,
\end{eqnarray*}
where
\begin{eqnarray}
\nabla_x\wedge p_b^{int}(R;x)&=&\left\{\begin{array}{ll}\mu_0q_j(x),&x\in\Omega_R^{int},
\\0,&x\in\Omega_R^{ext},\end{array}\right.\label{eqn:curlpbint}
\\\nabla_x\wedge p_b^{ext}(R;x)&=&\left\{\begin{array}{ll}0,&x\in\Omega_R^{int},
\\\mu_0q_j(x),&x\in\Omega_R^{ext}.\end{array}\right.\label{eqn:curlpbext}
\end{eqnarray}
In other words, $p_b^{int}(R;\cdot)$ denotes the part of the magnetic field that is due to source currents in the interior of the satellite's orbit $\Omega_R$, and $p_b^{ext}(R;\cdot)$ the part due to source currents in the exterior. A more detailed description can be found, e.g., in \cite{back96} and \cite{may06}. Since $p_b^{int}(R;\cdot)$ is still divergence-free, equation (\ref {eqn:curlpbint}) implies that a harmonic potential $U^{int}(R;\cdot):\Omega_R^{ext}\to\mathbb{R}$ exists, such that $p_b^{int}(R;x)=\nabla_xU^{int}(R;x)$, for $x\in\Omega_R^{ext}$. The potential $U^{int}$ can be expanded with respect to the outer harmonics $H_{n,k}^{ext}$, and the application of the gradient in combination with (\ref{eqn:hext}) then implies that $p_b^{int}(R;\cdot)$ can be expanded  in $\Omega_R^{ext}$ with respect to $\tilde{y}_{n,k}^{(1)}$. Analogously, $p_b^{ext}(R;\cdot)$ relates to an expansion in $\Omega_R^{int}$ with respect to $\tilde{y}_{n,k}^{(2)}$. The remaining toroidal part $q_b$ can be interpreted as the part induced by poloidal source currents $p_j$ crossing the sphere $\Omega_R$, and corresponds to the vector spherical harmonics $\tilde{y}_{n,k}^{(3)}$. To sum up, additionally observing that the poloidal and toroidal fields are continuous up to $\Omega_R$, we find 
\begin{eqnarray}\label{eqn:intextpol2}
b(x)=p_b^{int}(R;x)+p_b^{ext}(R;x)+q_b(x),\quad x\in\Omega_R,
\end{eqnarray}
with $\tilde{y}_{n,k}^{(i)}$, $i=1,2,3$, $n\in\mathbb{N}_{0_i}$, $k=1,\ldots,2n+1$, being the appropriate basis system for this split-up. Finally, the definition of the vector spherical harmonics in (\ref{eqn:ytilde}) implies that Theorem \ref{thm:althelmholtz} yields the exact same decomposition as (\ref{eqn:intextpol2}). More precisely, the $\tilde{o}^{(1)}$-part denotes the contribution due to sources in $\Omega_R^{int}$, the $\tilde{o}^{(2)}$-part the contribution due to sources in $\Omega_R^{ext}$, and the $\tilde{o}^{(3)}$-part the contribution due to source currents crossing the sphere $\Omega_R$. 

\subsection{Multiscale Representation}

Now, we turn to the actual multiscale representation. We discretize the regularized kernels from Theorem \ref{lem:althelmapprox}, by choosing parameters $\rho=2^{-J}$, for $J\in\mathbb{N}_{0}$. The \emph{scaling kernels (of scale $J$)} are then defined by 
\begin{eqnarray*}
\boldsymbol{\Phi}_J^{int}(\xi,\eta)&=&\boldsymbol{\Phi}_{2^{-J}}^{(1)}(\xi,\eta),\quad\xi,\eta\in\Omega,
\\\boldsymbol{\Phi}_J^{ext}(\xi,\eta)&=&\boldsymbol{\Phi}_{2^{-J}}^{(2)}(\xi,\eta),\quad\xi,\eta\in\Omega,
\\\boldsymbol{\Phi}_J^{q}(\xi,\eta)&=&\boldsymbol{\Phi}_{2^{-J}}^{(3)}(\xi,\eta),\quad\xi,\eta\in\Omega.
\end{eqnarray*}

\begin{figure}
\begin{center}
\scalebox{0.71}{\includegraphics{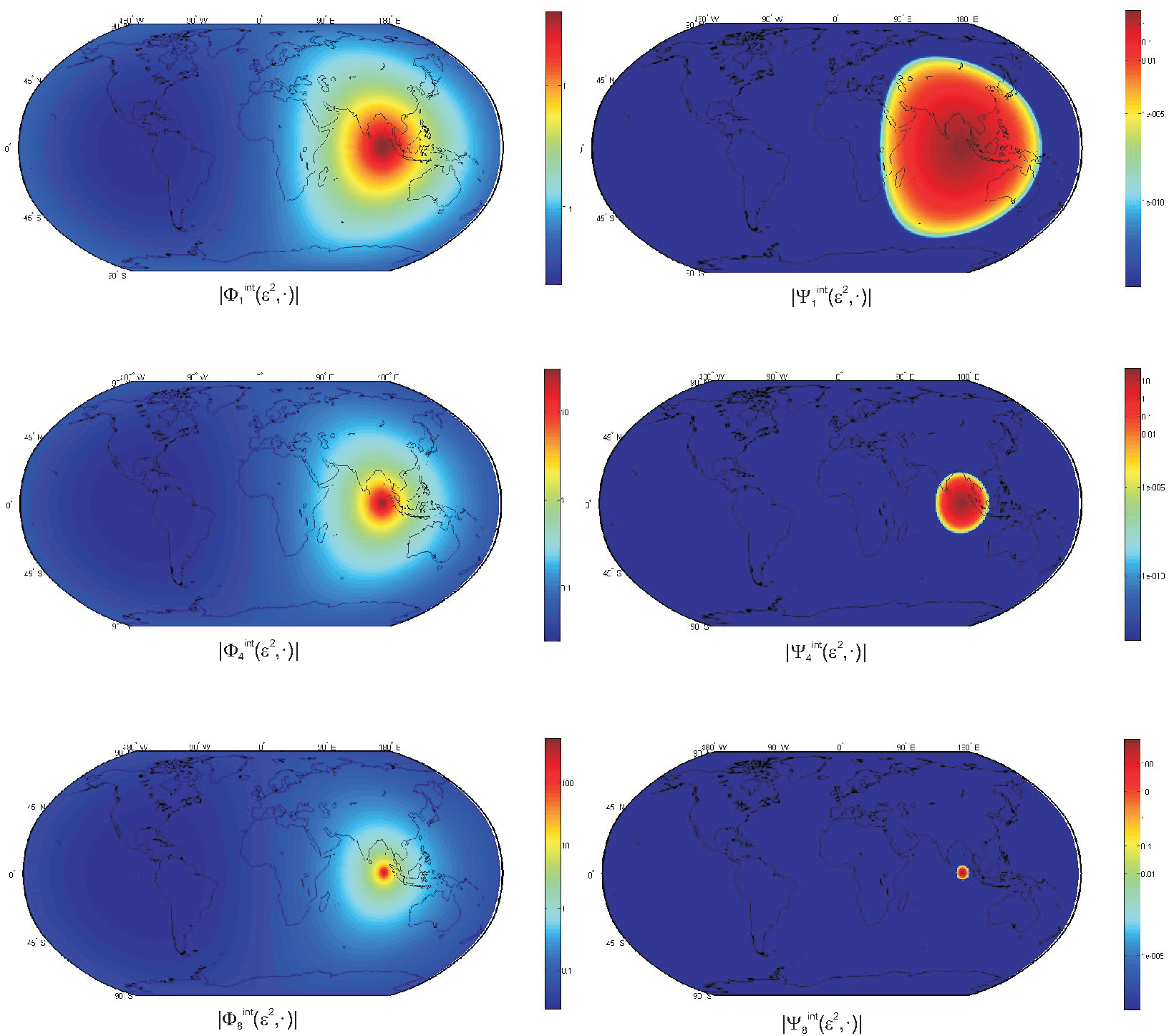}}
\end{center}
\caption{Absolute value of the scaling kernels (left) and wavelet kernels (right) for the interior poloidal contribution at scales $J=1,4,8$ (note that the color scaling is logarithmic).}\label{fig:kernels1}
\end{figure}

These kernels still have global support. The announced locally supported wavelets are obtained by taking the difference of two such scaling kernels. A \emph{wavelet kernel (of scale $J$)} denotes one of the following kernels
\begin{eqnarray*}
\boldsymbol{\Psi}_J^{int}(\xi,\eta)&=&\boldsymbol{\Phi}_{J+1}^{int}(\xi,\eta)-\boldsymbol{\Phi}_J^{int}(\xi,\eta),\quad\xi,\eta\in\Omega,
\\\boldsymbol{\Psi}_J^{ext}(\xi,\eta)&=&\boldsymbol{\Phi}_{J+1}^{ext}(\xi,\eta)-\boldsymbol{\Phi}_J^{ext}(\xi,\eta),\quad\xi,\eta\in\Omega,
\\\boldsymbol{\Psi}_J^{q}(\xi,\eta)&=&\boldsymbol{\Phi}_{J+1}^{q}(\xi,\eta)-\boldsymbol{\Phi}_J^{q}(\xi,\eta),\quad\xi,\eta\in\Omega.
\end{eqnarray*}
Due to the regularization of the Green function and the single layer kernel, these wavelets clearly have local support in a spherical cap of radius $2^{-J}$. More precisely, we find that $\textnormal{supp}\big( \boldsymbol{\Psi}_J^{i}(\xi,\cdot)\big)\subset\{\eta\in\Omega|1-\xi\cdot\eta<2^{-J}\}$, for $i\in\{int,ext,q\}$ and $\xi\in\Omega$. An illustration of the kernels is given in Figure \ref{fig:kernels1}.

Each scaling kernel generates a scaling transform. These are given by
\begin{eqnarray}
P_J^{i}b(x)&=&\int_\Omega \boldsymbol{\Phi}^{i}_J(\xi,\eta) b(R\eta)d\omega(\eta),\quad x=R\xi\in\Omega_R.\label{eqn:scaltransformint}
\end{eqnarray}
The corresponding wavelet transforms read
\begin{eqnarray}
R_J^{i}b(x)&=&\int_\Omega \boldsymbol{\Psi}^{i}_J(\xi,\eta) b(R\eta)d\omega(\eta),\quad x=R\xi\in\Omega_R.\label{eqn:wavelettransformint}
\end{eqnarray}
The idea of the multiscale approach is to resolve the modeled quantities at different spatial resolutions. Therefore, the scaling kernels are only used to provide a trend approximation of the coarse features at some small initial scale $J_0$. The spatially stronger localized features are subsequently added by the wavelet transforms. This is reflected in the following relations,
\begin{eqnarray}
P_J^{i}b(x)&=&P_{J_0}^{i}b(x)+\sum_{j=J_0}^{J-1}R_j^{i}b(x),\quad x=R\xi\in\Omega_R,
\end{eqnarray}
for $i\in\{int,ext,q\}$. Different from the multiresolution constructed, e.g., in \cite{may06} and \cite{maymai}, the scale spaces $V_J^i=\{P_J^ib|b\in c^{(2)}(\Omega_R)\}$ in our approach are not necessarily nested in the sense $V_J^i\subset V_{J+1}^i$. The advantage here is the local support of the wavelet kernels, which implies that the evaluation of the wavelet transforms $R_j^{i}b(x)$, for $x=R\xi\in\Omega_R$, only requires data in a spherical cap around $\xi\in\Omega$ with scale-dependent spherical radius $2^{-j}$. Thus, regions with a higher data density can be resolved up to higher scales, i.e., up to a higher spatial resolution, without suffering errors from the lower data densities in surrounding areas. The general concept is illustrated by the following tree algorithm
\begin{center}
\unitlength0.9mm
\begin{picture}(150,20)
\put(0,0){\makebox(0,0)[c]{$P^i_{J_0}F$}}
\put(20,0){\makebox(0,0)[c]{$+$}}
\put(40,0){\makebox(0,0)[c]{$P^i_{J_0+1}F$}}
\put(62,0){\makebox(0,0)[c]{$+$}}
\put(82,0){\makebox(0,0)[c]{$P^i_{J_0+2}F$}}
\put(104,0){\makebox(0,0)[c]{$+$}}
\put(123,0){\makebox(0,0)[c]{$\ldots$}}
\put(0,15){\makebox(0,0)[c]{$R^i_{J_0}F$}}
\put(40,15){\makebox(0,0)[c]{$R^i_{J_0+1}F$}}
\put(82,15){\makebox(0,0)[c]{$R^i_{J_0+2}F$}}
\put(5,0){\vector(1,0){12}}
\put(47,0){\vector(1,0){12}}
\put(89,0){\vector(1,0){12}}
\put(23,0){\vector(1,0){10}}
\put(65,0){\vector(1,0){10}}
\put(107,0){\vector(1,0){10}}
\put(5,13){\vector(1,-1){12}}
\put(47,13){\vector(1,-1){12}}
\put(89,13){\vector(1,-1){12}}
\put(128,0){\vector(1,0){10}}
\put(146,0){\makebox(0,0)[c]{$P^i_{J_{max}}F$.}}
\end{picture}\\[0ex]\textbf{}
\end{center}
The maximal scale $J=J_{max}$ at which $P_Jb(x)$ can be evaluated is determined by the amount of data points in the vicinity of $x=R\xi\in\Omega_R$. Sufficiently many data points in the support of $\boldsymbol{\Psi}^{i}_{J_{max}-1}(\xi,\cdot)$ are required to guarantee a numerical meaningful evaluation of the integral in the wavelet transform $R_{J_{max}-1}^{i}b(x)$.  

To conclude this subsection, we summarize the results in the following theorem, which is mainly a reformulation of Theorem \ref{lem:althelmapprox} in terms of the above described multiscale setting for the separation of the magnetic field with respect to its sources.

\begin{thm}\label{thm:intpart}
Let $b$ be of class $c^{(2)}(\mathbb{R}^3)$, and $j$ of class $c^{(1)}(\mathbb{R}^3)$, satisfying the pre-Maxwell equations
\begin{eqnarray*}
\nabla_x\wedge b(x)&=&\mu_0j(x),\quad x\in\mathbb{R}^3,
\\\nabla_x\cdot b(x)&=&0,\quad x\in\mathbb{R}^3.
\end{eqnarray*}
If $P_J^{int}$, $P_J^{ext}$, $P_J^{q}$, $R_J^{int}$, $R_J^{ext}$, $R_J^{q}$ are defined as in (\ref{eqn:scaltransformint}), (\ref{eqn:wavelettransformint}), then 
\begin{eqnarray*}
b(x)&=&p_b^{int}(R;x)+p_b^{ext}(R;x)+q_b(x),\quad x\in\Omega_R,
\end{eqnarray*}
for a fixed $R>0$, with
\begin{eqnarray*}
p_b^{int}(R;x)&=&P_{J_0}^{int}b(x)+\sum_{j=J_0}^{\infty}R_j^{int}b(x),\quad x\in\Omega_R,
\\p_b^{ext}(R;x)&=&P_{J_0}^{ext}b(x)+\sum_{j=J_0}^{\infty}R_j^{ext}b(x),\quad x\in\Omega_R,
\\q_b(x)&=&P_{J_0}^{q}b(x)+\sum_{j=J_0}^{\infty}R_j^{q}b(x),\quad x\in\Omega_R.
\end{eqnarray*}\end{thm}

\subsection{Crustal Field Modeling from CHAMP Data}

\begin{figure}[b]
\begin{center}
\scalebox{0.73}{\includegraphics{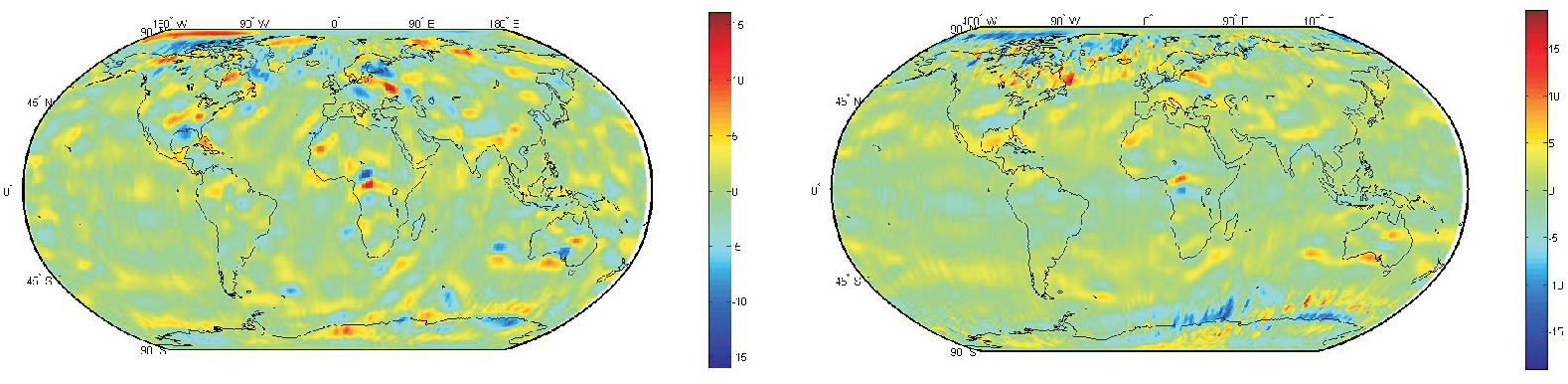}}
\end{center}
\caption{The radial component (left) and the south-north component (right) of the input magnetic field (in nT), averaged to a $180\times 180$ equiangular grid.}\label{fig:binput}
\end{figure}

In this subsection, we apply the above derived multiscale approach to a set of CHAMP satellite measurements. The used data set is similar to the one used in \cite{may06} and \cite{maymai} and has been collected between June 2001 and December 2001. It has been pre-processed at the GFZ Potsdam by Stefan Maus to fit the purpose of crustal field modeling (see, e.g., \cite{maus06} for a detailed description). Due to the almost spherical orbit of the CHAMP satellite, we can assume all data to be given on a sphere of radius $R_E+450$km, where $R_E=6371.2$km denotes the mean Earth radius. For the discretization of the integrals appearing in the multiscale representation, we use the integration rule described in \cite{drihea}, which requires an equiangular data grid (and reflects the data situation of satellite measurements). In this example, we use a grid with $180$ grid points in latitudinal as well as longitudinal direction. Centered around each grid point, we select a spherical rectangle with a diameter of $2.5^\circ$ in latitude and longitude and average all measurements in the cell using a M-estimation with Huber's weight function (see, e.g., \cite{hogg}). The resulting input data set is shown in Figure \ref{fig:binput}.

\begin{figure}
\begin{center}
\scalebox{0.72}{\includegraphics{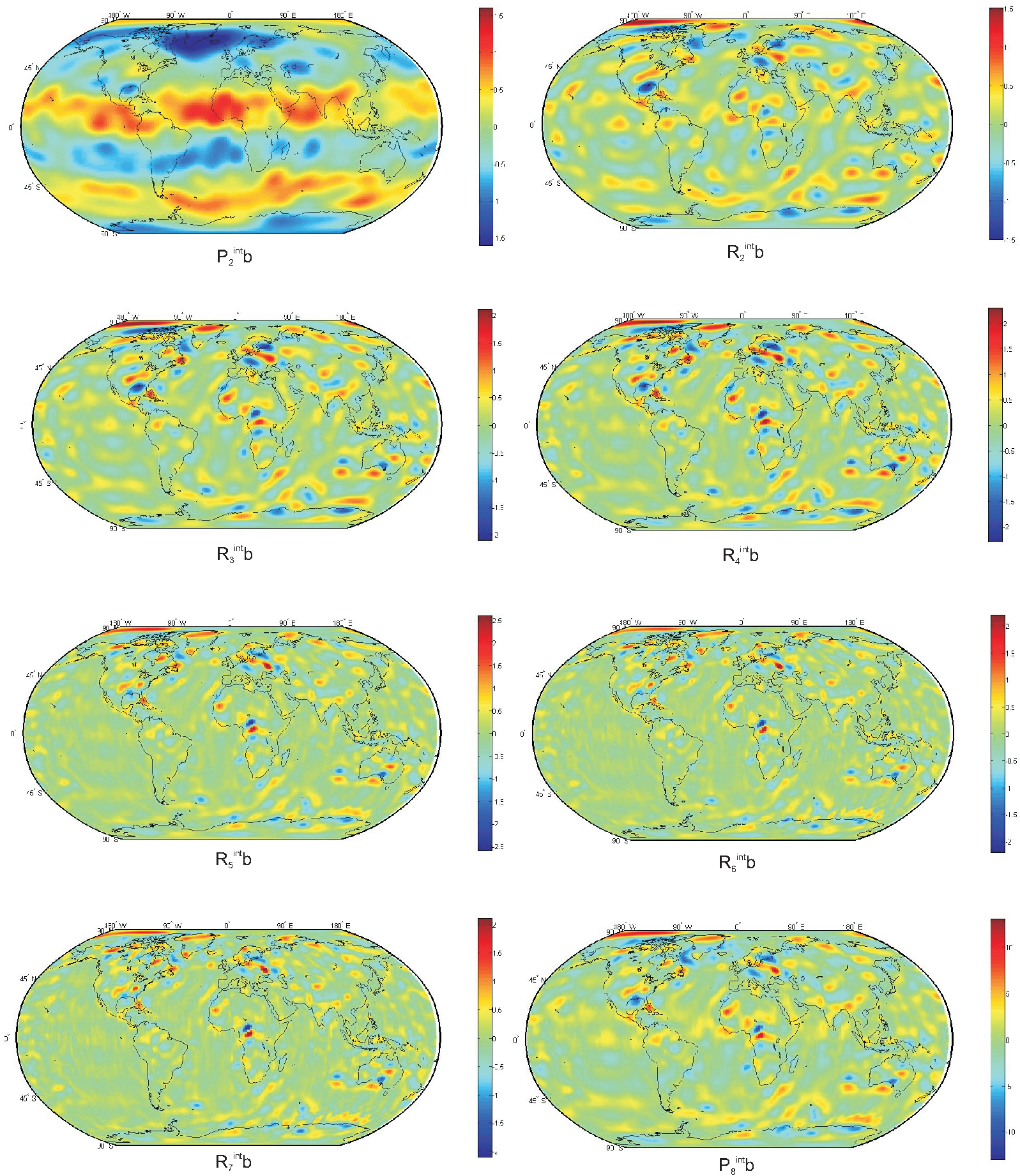}}
\end{center}
\caption{Approximation of the internal poloidal magnetic field (in nT), at initial scale $J=2$ (top left) and at scale $J=8$ (bottom right); only the radial component is shown. The remaining figures show the intermediate wavelet contributions from scale $J=2$ to scale $J=7$ (top right to bottom left).}\label{fig:brmultiscal}
\end{figure}

\begin{figure}
\begin{center}
\scalebox{0.72}{\includegraphics{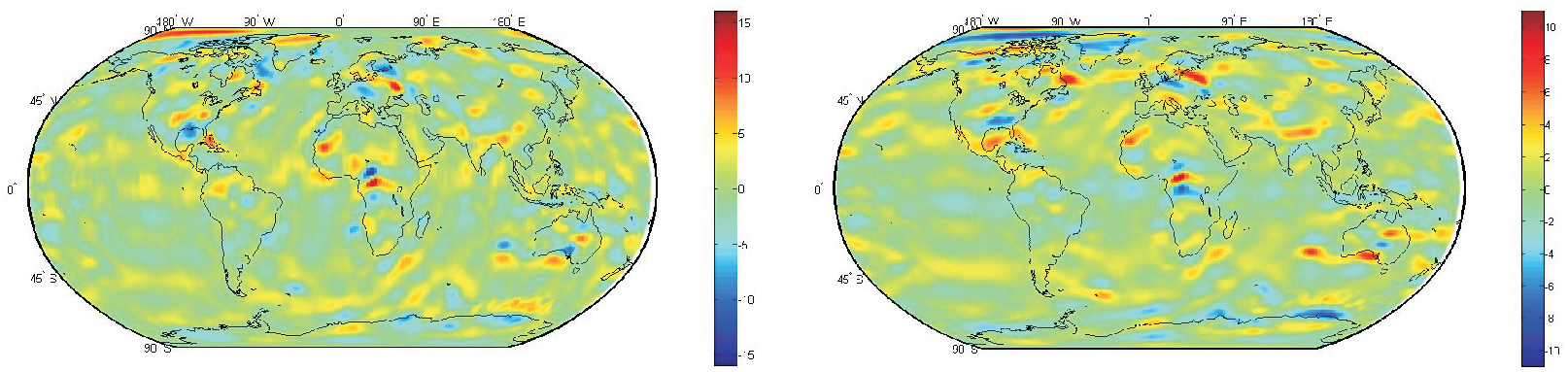}}
\end{center}
\caption{Approximation of the internal poloidal magnetic field (in nT), at scale $J=9$; radial component (left) and south-north component (right).}\label{fig:bs9}
\end{figure}

\begin{figure}
\begin{center}
\scalebox{0.72}{\includegraphics{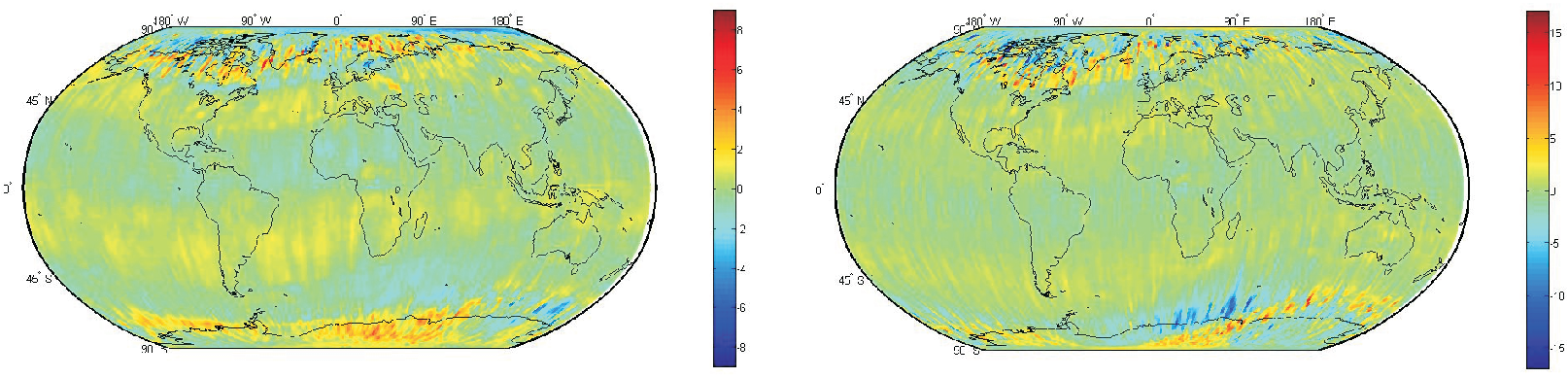}}
\end{center}
\caption{Difference between the input data set and the approximation of the internal poloidal magnetic field (in nT), at scale $J=9$; radial component (left) and south-north component (right).}\label{fig:bs9diff}
\end{figure}

The results obtained from the multiscale representation are illustrated in Figures \ref{fig:brmultiscal}--\ref{fig:bs9diff}. For the sake of brevity, we only indicate the radial and the south-north component of the magnetic field, and in Figure \ref{fig:brmultiscal} only the radial component. Figure \ref{fig:brmultiscal} also illustrates best the different spatial resolutions of the multiscale representation.   
The initial trend approximation at scale $J_0=2$ only resolves very coarse features, while the subsequent wavelet transforms resolve more and more localized features, such that the scales $J=5,6,7$ mainly focus on the strongest crustal field anomalies over Central Africa and Eastern Europe, as well as North America and Australia. In oceanic regions, there is hardly any contribution at these scales, indicating that there the crustal field is of a rather coarse nature (i.e., of large wavelength when arguing in frequency domain). Furthermore, one finds that the structure of the wavelet contributions hardly changes for scales higher than $J=5$. This might be an indicator of the general spatial extend of the anomalies of the crustal field signal at satellite altitude (when comparing the resolved features with the size of the support of the wavelet kernels). Figure \ref{fig:bs9} shows the final approximation $P_{J_{max}}^{int}b$ of the internal magnetic field contributions at the highest scale $J_{max}=9$. The difference to the input data set is indicated in Figure \ref{fig:bs9diff}. It actually illustrates the performance of the separation with respect to the sources. One can recognize strong polar fields that are clearly not due to the Earth's crustal field and are probably induced by polar ionospheric current systems. Furthermore, one finds bands of positive and negative field strength oriented parallel to the dipole equator. This is a typical signature of magnetospheric ring currents. 

Thus, the multiscale approach of this paper can be used to improve pre-processed crustal magnetic field data. The decomposition with respect to $\tilde{o}^{(i)}$ is actually able to filter out contributions originating outside the satellite's orbit or at satellite altitude, while the wavelet transformations give a clearer impression of the local features of the crustal field.


\end{document}